\documentclass[12pt]{article}
\usepackage[utf8]{inputenc}
\usepackage{geometry}
\usepackage{epsfig}
\usepackage{oldgerm}
\usepackage{amsmath}
\usepackage{amssymb}
\usepackage{amsthm}
\usepackage{cite}
\usepackage{enumerate}
\usepackage{amstext,color}
\usepackage{fullpage}
\usepackage{graphicx,graphics}
\usepackage{graphicx}
\usepackage{algpseudocode} 
\usepackage{algorithm} 
\usepackage{float}
\usepackage{subcaption}

\usepackage[caption = false]{subfig}
\DeclareMathOperator*{\esssup}{ess\,sup}
\newtheorem{theorem}{Theorem}[section]
\newtheorem{definition}{Definition}[section]
\newtheorem{lemma}{Lemma}[section]

\newtheorem{ass}{Assumption}

\newtheorem{remark}{Remark}[section]
\newcommand{\x}[0]{\boldsymbol{x}}
\usepackage{ragged2e}
\justifying
\begin{document}
	\begin{center}
		{{\bf \large {\rm {\bf Space-Time Isogeometric Method for a Nonlocal Parabolic Problem }}}}
	\end{center}
	\begin{center}
		
		{\textmd {\bf Sudhakar Chaudhary}}\footnote{\it Department of Basic Sciences, Institute of Infrastructure, Technology, Research and Management, Gujarat, India, (dr.sudhakarchaudhary@gmail.com)}
		{\textmd {{\bf Shreya Chauhan}}}\footnote{\it  Department of Basic Sciences, Institute of Infrastructure, Technology, Research and Management, Gujarat, India, (shreyachauhan898@gmail.com)}
    {\textmd {{\bf Monica Montardini* }}}\footnote{\it *Corresponding Author, Dipartimento di Matematica, Università degli Studi di Pavia, via Ferrata 5, Pavia, Italy, (monica.montardini@unipv.it)}
	\end{center}
\begin{abstract}In the present work, we focus on the space-time isogeometric discretization of a parabolic problem with a nonlocal diffusion coefficient. 
The existence and uniqueness of the solution for the continuous space-time variational formulation are proven. We prove the existence of the discrete solution and also establish the a priori error estimate for the space-time isogeometric scheme. The non-linear system is linearized through Picard's method and a suitable preconditioner for the linearized system is provided.  Finally, to confirm the theoretical findings, the results of some numerical experiments are presented.
	\end{abstract}
	\begin{center}
	\end{center}
\section{Introduction}\label{sec1}
\par In this paper, the following nonlocal problem is studied: Find $u(\boldsymbol{x},t)$ such that,
\begin{equation}\label{eqn1.1}
		\begin{split}
			\displaystyle \partial_tu(\x,t) -a(l(u)(t))\Delta u(\x,t) &= f(\x,t)\quad \mbox{in } \Omega\times(0,T), \\
			u(\boldsymbol{x},t) &=0\quad\mbox{on } \partial\Omega\times(0,T),  \\
			u(\boldsymbol{x},0)&=0\quad\mbox{in } \Omega,
		\end{split}
	\end{equation}
	where $\partial_tu=\frac{\partial u}{\partial t},$  $\Omega$ is a bounded domain with Lipschitz boundary $\partial\Omega$ in $\mathbb{R}^d$, $d= 1,2,3$ and $f:\Omega\times(0,T)\rightarrow \mathbb{R}$ and $a:\mathbb{R}\rightarrow\mathbb R^+$ are given functions. The term $l(u)$ is chosen as  in \cite{Chipot1999}:
	\begin{equation*}
		l(u)=\int_{\Omega}u\ d\Omega .
	\end{equation*}
\par In recent years, the study of nonlocal problems has received great attention due to its usefulness in real-world applications (see \cite{app2,app3,app4}). Problem \eqref{eqn1.1} is used to model the population density of bacteria in a given vessel (the domain $\Omega$) subject to spreading or propagation of heat \cite{Chipot1999}. The well-posedness of the solution to problem \eqref{eqn1.1} and its stationary counter part was studied in \cite{CHIPOT19974619}. 
In \cite{Chipot1999}, authors established existence-uniqueness of a weak solution to problem \eqref{eqn1.1} using  Faedo–Galerkin approximations and compactness arguments, and they  described also the asymptotic behaviour of the solution. In \cite{dirichipot}, the well-posedness of problem  \eqref{eqn1.1} was studied with the choice $l(u)= \int_{\Omega}\vert\nabla u\vert^2\ d\Omega $.
\par To obtain a numerical solution of problems like \eqref{eqn1.1}, the usual time-stepping methods combined with various spatial discretization techniques such as the Finite Element Method (FEM) \cite{chaudhary2017finite, ROBALO20145609, ALMEIDA2016146} and Virtual Element Method (VEM) \cite{Adak2020} are used. This approach is sequential in time, complicating the parallelization of the algorithms and affecting the simultaneous adaptive refinement   in space and time. As an alternative to time-stepping methods, an overall space-time discretized problem can be considered. We refer to \cite{langer2019space} for more details about space-time finite element methods. In \cite{Steinbach+2015+551+566}, a completely unstructured space-time finite element method was analyzed to solve a linear parabolic problem. The author in \cite{moore2018stable} proposed a stable space-time finite element method with time-upwind test functions for the numerical solution of linear evaluation problems in moving domain. Recently in \cite{p-lap,plap2}, the space-time finite element scheme was analyzed for the nonlinear parabolic $p$-Laplace problem.
\par Pioneered by Hughes et.al. in \cite{HUGHES20054135},  Isogeometric Analysis (IgA) is developed on the idea to use the same functions (typically B-splines or Non-Uniform Rational B-splines (NURBS)) for the parameterization of the domain and for the approximation of the solution. The main aim behind IgA is to increase the interoperability between computer-aided design industry and numerical simulations.
\par A full space-time approach with isogeometric method was introduced by Langer et. al. in \cite{stiga1}.  In \cite{LOLI20202586}, the authors developed a suitable preconditioner and an efficient solver for a Galerkin space-time isogeometric discretization of a linear parabolic problem. An efficient solver for a least-squares discretization of a parabolic problem is provided in \cite{leastsquare}. A space-time isogeometric method for linear and non-linear electrodynamics problems was investigated in \cite{Saad2021}. In \cite{stiganonlinear1}, a stabilized space-time isogeometric method has been designed in the context of cardiac electrophysiology.
\par Theoretically, the space-time method is well studied for linear time-dependent problems. However, the literature on the space-time method for non-linear time-dependent problems is very sparse. In particular, convergence analysis of space-time methods for nonlinear problems is not much available. In our work, we have attempted to use the space-time isogeometric method to solve a nonlinear nonlocal parabolic problem and establish the convergence results theoretically. Another reason behind considering space time IgA for nonlocal problem is that the classical FEM uses low degree polynomials and the geometry of the problem is represented approximately. As an alternative, IgA allows exact representation of the geometry, and the use of smooth splines can enhance the accuracy of the numerical solution. Moreover, compared with time-stepping methods, the space-time method provides additional benefits such as parallelization in time and high-order accuracy both in space and time. To handle the nonlinear term, we use Picard's iterative scheme. Each linearized iteration is solved using  GMRES method, coupled with a preconditioner that is a variant of \cite{LOLI20202586}. To the best of our knowledge, this is the first work on a space-time isogeometric method for nonlocal parabolic problems.



\par The remainder of the paper is outlined as follows: Some preliminaries on Isogeometric Analysis and function spaces are illustrated in Section \ref{secPrel}. In Section \ref{sec2}, we write the space-time weak formulation of problem \eqref{eqn1.1} and we give some results on the existence and uniqueness of the weak solution.  The space-time isogeometric discretization of \eqref{eqn1.1} is given in Section \ref{sec4}. In this section, we also provide a suitable solver for the linearized problem. The error estimate of the numerical scheme is presented in Section \ref{sec5}. Finally, we provide some numerical experiments to verify the theoretical findings in Section \ref{sec7}.

\section{Preliminaries}
\label{secPrel}
We begin by defining some function spaces and then we introduce the isogeometric spaces for the numerical approximation of \eqref{eqn1.1}. We also recall some results from functional analysis.
\subsection{Function spaces}
First, we define some function spaces which will be used in our analysis. Let $L^2(\Omega)$ be the Hilbert space of all square integrable functions with norm $\vert\vert v\vert\vert_{L^2(\Omega)}=\left(\int_{\Omega}\vert v(\boldsymbol{x})\vert ^2\ d\Omega \right)^{\frac{1}{2}}$ and inner product $(w,v)=\int_{\Omega}w(\boldsymbol{x})v(\boldsymbol{x}) \ d\Omega$  $\forall w,v\in L^2(\Omega)$. Let $H^r(\Omega)$ denote the Sobolev space on $\Omega$ for $r\in \mathbb N$, with norm 
  \begin{equation*}
      \left\| v\right\|_r=\left(\sum_{0\leq \vert\alpha\vert\leq r}\left\lVert\frac{\partial^{\alpha}v}{\partial \boldsymbol{x}^{\alpha}}\right\rVert^2\right)^{\frac{1}{2}}.
  \end{equation*}
  For $1\leq p\leq \infty$, by $L^p(0,T;H^r(\Omega))$, we denote the space of the vector valued functions $v:(0,T)\rightarrow H^r(\Omega)$, with norm $\vert\vert v\vert\vert_{L^p(0,T;H^r(\Omega))}=\left(\int_{0}^{T}\vert\vert v(t)\vert\vert_r^p\ dt\right)^{\frac{1}{p}}$ for $1\leq p<\infty$ and $\vert\vert v\vert\vert_{L^{\infty}(0,T;H^r(\Omega))}=\displaystyle\esssup_{0\leq t\leq T}\vert\vert v(t)\vert\vert_r$ for $p=\infty$. To derive a space-time weak formulation of \eqref{eqn1.1}, we consider the following space-time function spaces: $V=\{v\in  L^2(0,T;H_0^1(\Omega))\ | \ \partial_tv\in L^2(0,T;H^{-1}(\Omega))\ \text{ and }\ v(\boldsymbol{x},0)=0\}$ and $W=L^2(0,T;H_0^1(\Omega))$ with the corresponding norms given by
	\begin{equation*}
		||v||^2_V:=||\partial_tv||^2_{L^2(0,T;H^{-1}(\Omega))}+||v||^2_{L^2(0,T;H_0^1(\Omega))}
	\end{equation*}
	and
	 \begin{equation*}
		||v||^2_W:=||v||^2_{L^2(0,T;H_0^1(\Omega))}.
	\end{equation*} 
	Here, $H^{-1}(\Omega)$ denotes the dual space of $H_0^1(\Omega)$. We note that the dual of the space $W$ is $W^*=L^2(0,T;H^{-1}(\Omega))$. We also consider the space $Z=L^2(0,T;L^2(\Omega))$ with norm 
 \begin{equation*}
		||v||^2_Z:=||v||^2_{L^2(0,T;L^2(\Omega))}.
	\end{equation*} 
Throughout the paper, $C$ and $C'$ are generic constants, and $C_p$ denotes the Poincar\'e constant.
\subsection{Isogeometric spaces}\label{sec3}
	Some fundamental aspects of B-spline spaces and isogeometric concepts are given in the present subsection. For further discussion on B-splines, we refer to the monograph \cite{NURBSBook}, while for an insight on mathematical aspects of IgA, we refer to \cite{BEIRAODAVEIGA20121}.\\
	Let $n$ and $q$ be two non-negative integers. A knot vector in $[0,1]$ is a  set of non-decreasing real numbers  
	\begin{equation*}
		\Xi=\{0=\zeta_1\leq\cdots\leq\zeta_{n+q+1}=1\}.
	\end{equation*}
	 We will consider only open knot vectors, i.e. we set  $\zeta_1=\cdots=\zeta_{q+1}$ and $\zeta_{n+1}=\cdots=\zeta_{n+q+1}$. The $j^{th}$ B-spline basis function is defined   using Cox-de Boor recursion formula as:
	\begin{equation*}
		\widehat{b}_{j,0}(\zeta)=\begin{cases}
			1 & \mbox{if } \zeta_j\leq \zeta<\zeta_{j+1}, \\
			0 & \mbox{otherwise}
		\end{cases}
	\end{equation*}
 and
	\begin{equation*}
		\displaystyle\widehat{b}_{j,q}(\zeta)=\frac{\zeta-\zeta_j}{\zeta_{j+q}-\zeta_j}\widehat{b}_{j,q-1}(\zeta)+\frac{\zeta_{j+q+1}-\zeta}{\zeta_{j+q+1}-\zeta_{j+1}}\widehat{b}_{j+1,q-1}(\zeta)\ \text{for }q\geq 1,
	\end{equation*}
	where we take $0/0=0$.  By $h$, we denote the meshsize which is defined as $h:=\text{max}\{|\zeta_{j+1}-\zeta_j|: j=1,\cdots,n+q\}$. Then, we define the spline space of functions in one variable as
	\begin{equation*}
		\widehat{\mathcal{S}}_h:=\text{span}\left\{\widehat{b}_{j,q} \mid  j=1,\cdots,n\right\}.
	\end{equation*}
	The continuity of the B-splines at knots depends on the multiplicity of the knots (see \cite{NURBSBook}). In this work we only consider maximum regularity splines, as the regularity of splines improves the isogeometric approximation per degree of freedom \cite{da2014mathematical}. \\
Multivariate B-splines can be defined by taking the tensor product of the univariate B-splines. Let $n_r,q_r$ for $1\leq r\leq d$ and $n_t,q_t$ be positive integers. Then we consider $d+1$ open knot vectors $\Xi_r:=\{\zeta_{r,1},\cdots,\zeta_{r,n_r+q_r+1}\}$ for $1\leq r\leq d$ and $\Xi_t:=\{\zeta_{t,1},\cdots,\zeta_{t,n_t+q_t+1}\}$. We denote by $h_t$, the meshsize associated with $\Xi_t$ and $h_r$ the meshsize associated with $\Xi_r$ for $1\leq r\leq d$. Let $h_s:=\text{max}\{h_r,1\leq r\leq d\}$ be the maximum of all meshsizes in each spatial direction. We define the degree vector by $\mathbf{q}=(\mathbf{q}_s,q_t)$ with $\mathbf{q}_s=(q_1,\cdots,q_d)$. In the present work, for simplicity, we consider the same degree in all the spatial directions and we set, with a little abuse of notations, $q_s:=q_1=\cdots=q_d$. 
	We  consider quasi-uniform knot vectors, i.e. we suppose that  there exists $0<\alpha<1$, such that $\alpha h_s\leq \zeta_{r,j+1}-\zeta_{r,j}\leq h_s$  for every knot span $(\zeta_{r,j},\zeta_{r,j+1})$ with $\zeta_{r,j}\neq\zeta_{r,j+1}$ of $\Xi_r$ and for $1\leq r\leq d$ and every knot span $(\zeta_{t,j},\zeta_{t,j+1})$ with $\zeta_{t,j}\neq\zeta_{t,j+1}$ of $\Xi_t$ satisfies $\alpha h_t\leq \zeta_{t,j+1}-\zeta_{t,j}\leq h_t$.
	Multivariate B-splines   are defined as:
	\begin{equation*}
		\widehat{B}_{\mathbf{j},\mathbf{q}}(\mathbf{\zeta},\tau):=\widehat{B}_{\mathbf{j}_s,\mathbf{q}_s}(\mathbf{\zeta})\widehat{b}_{j_t,q_t}(\tau),
	\end{equation*}
	where
	\begin{equation*}
		\widehat{B}_{\mathbf{j}_s,\mathbf{q}_s}(\mathbf{\zeta})=\widehat{b}_{j_1,q_s}(\zeta_1)\cdots\widehat{b}_{j_d,q_s}(\zeta_d),
	\end{equation*}
	with $\mathbf{j}_s=(j_1,\cdots,j_d)$, $\mathbf{j}=(\mathbf{j}_s,j_t)$ and  $\mathbf{\zeta}=(\zeta_1,\cdots,\zeta_d)$.\\
	We define the corresponding spline function space as
	\begin{equation*}
		\widehat{\mathcal{S}}_{h,\mathbf{q}}:=\text{span}\left\{\widehat{B}_{\mathbf{j},\mathbf{q}} \mid  j_k=1,\dots , n_k \text{ and } k=1,\dots, d;  j_t=1,\dots , n_t\right\},
	\end{equation*}
	where $h=\text{max}\{h_s,h_t\}$.\\
	We can write $ \widehat{\mathcal{S}}_{h,\mathbf{q}}$ as a tensor product of the spatial spline space and the spline space for the time variable as
	\begin{equation*}
		\widehat{\mathcal{S}}_{h,\mathbf{q}}= \widehat{\mathcal{S}}_{h_s,\mathbf{q_s}}\otimes \widehat{\mathcal{S}}_{h_t,{q_t}},
	\end{equation*}
	where
	\begin{equation*}
		\widehat{\mathcal{S}}_{h_s,\mathbf{q_s}}:=\text{span}\left\{\widehat{B}_{\mathbf{j}_s,q_s} \mid  j_k=1,\dots , n_k \text{ and } k=1,\dots, d\right\}
	\end{equation*}
 and
 \begin{equation*}
    \widehat{\mathcal{S}}_{h_t,{q_t}}:=\text{span}\{\widehat{b}_{{j}_t,q_t}  \mid  j_t=1,\dots,n_t \}.
 \end{equation*}
 We need the following assumption.
\begin{ass}
We assume that $q_t$, $q_s\geq 1$  and that $\widehat{\mathcal{S}}_{h_s,\mathbf{q_s}}\subset C^0(\widehat{\Omega})$ and $ \widehat{\mathcal{S}}_{h_t,{q_t}}\subset C^0((0,1))$.
\end{ass}
 
	Now, we define isogeometric function spaces on the space-time cylinder $\Omega\times(0,T)$. The following assumption is considered in our analysis.

	\begin{ass}  Let $\widehat{\mathbf{F}}:\widehat{\Omega}\rightarrow\Omega$ be the   parametrization of  the spatial domain $\Omega$.  We suppose that $\widehat{\mathbf{F}}\in [\widehat{\mathcal{S}}_{h_s,\mathbf{q_s}}]^d$ and that the piece-wise derivatives of any order of $\widehat{\mathbf{F}}^{-1}$ are bounded.\end{ass}

The space-time cylinder is parameterized by the map $\widehat{\mathbf{G}}:\widehat{\Omega}\times(0,1)\rightarrow\Omega\times(0,T)$, $\widehat{\mathbf{G}}(\mathbf{\zeta},\tau):=(\widehat{\mathbf{F}}(\mathbf{\zeta}),T\tau)$.\\
	To incorporate the boundary and initial conditions, we consider the following spline space in parametric coordinates  
	\begin{equation*}
		\widehat{V}_h:=\{\widehat{u}_h\in\ \widehat{\mathcal{S}}_{h,\mathbf{q}}\ | \ \widehat{u}_h=0\text{ on } \partial\widehat{\Omega}\times(0,1)\text{ and }\widehat{u}_h=0\text{ in }\widehat{\Omega}\times\{0\}\}.
	\end{equation*}
	Due to the tensor product structure, we have $\widehat{V}_h=\widehat{V}_{s,h_s}\otimes\widehat{V}_{t,h_t}$,
	where
	\begin{equation*}
		\widehat{V}_{s,h_s}=\{\widehat{v}_h\in \widehat{\mathcal{S}}_{h_s,\mathbf{q_s}}\ | \ \widehat{v}_h=0\text{ on }\partial\widehat{\Omega}\}=\text{span}\{\widehat{B}_{j,\mathbf{q_s}}\ | \ j=1,\cdots,N_s\}
	\end{equation*}
	and
	\begin{equation*}
		\widehat{V}_{t,h_t}:=\{\widehat{v}_h\in \widehat{\mathcal{S}}_{h_t,{q_t}}\ | \ \widehat{v}_h(0)=0\}=\text{span}\{\widehat{b}_{j,q_t}\ | \ j=1,\cdots,N_t\},
	\end{equation*}
 where  $N_s=\prod_{l=1}^{d} n_{s,l}$, $n_{s,l}=n_l-2$ for $1\leq l\leq d$, $N_t=n_t-1$ and where we have used a colexicographical reordering of the degrees of freedom.
	Then
	\begin{equation*}
		\widehat{V}_h=\text{span}\{\widehat{B}_{j,\mathbf{q}}\ | \ j=1,\cdots,N_{dof}\},
	\end{equation*}
	where $N_{dof}=N_sN_t$.\\
	Finally, the isogeometric space is given by
	\begin{equation}
 \label{eq:isospace}
		V_h=\text{span}\{B_{j,\mathbf{q}}:=\widehat{B}_{j,\mathbf{q}}\circ\widehat
  {\mathbf{G}}^{-1}\ | \ j=1,\cdots,N_{dof}\}
	\end{equation}
 and we can write
 \begin{equation*}
     V_h=V_{s,h_s}\otimes V_{t,h_t},
\end{equation*}
where
\begin{equation*}
    V_{s,h_s}=\text{span}\{B_{j,\mathbf{q}_s}=\widehat{B}_{j,\mathbf{q}_s}\circ\widehat{\mathbf{F}}^{-1}\ | \ j=1,\cdots,N_s\}
\end{equation*}
and
\begin{equation*}
    V_{t,h_t}=\text{span}\{b_{j,q_t}=\widehat{b}_{j,q_t}(\cdot/T):j=1,\cdots,N_t\}.
\end{equation*}
Next, we report a result from \cite{LOLI20202586} regarding the approximation property of the isogeometric space $V_h$, which will be useful in proving the a priori error estimate.
\begin{theorem}\cite[Theorem 2]{LOLI20202586}\label{projerr}
Let $r$ be an integer such that $1< r\leq \min\{q_s,q_t\}+1$ and $u\in V\cap H^r(\Omega\times(0,T))$. Then there exists a projection $\mathcal{P}_h: V\cap H^r(\Omega\times(0,T))\rightarrow V_h$ such that 
\begin{equation}\label{eqprojerr}
    \vert\vert u-\mathcal{P}_h u\vert\vert_V\leq C (h_t^{r-1}+h_s^{r-1})\vert\vert u\vert\vert_{H^r(\Omega\times(0,T))},
\end{equation}
  where the constant $C$ is independent of the meshsizes $h_s$ and $h_t$.  
\end{theorem}
\subsection{Special types of operators}
In this subsection, we recall some definitions and results of functional analysis.
By $X$, we denote a reflexive Banach space and $X^{*}$ denotes its dual. The duality pairing between $X^*$ and $X$ is denoted by the symbol $\langle\cdot ,\cdot\rangle.$
\begin{definition}\cite[Section 2.2, Section 2.3]{lsimon}
	Let $A:X\rightarrow X^{*}$ be an operator. We say that
	\begin{itemize}
		\item $A$ is bounded if each bounded subset of $X$ is mapped to a bounded subset of $X^{*}$.
		\item $A$ is demicontinuous if for every sequence $(v_k)_k$ in $X$ converging to $v$ in $X$, $(A(v_k))_k$ converges weakly to $ A(v)$ in $X^{*}$.
		\item $A$ is monotone if $\langle A(v_1)-A(v_2),v_1-v_2\rangle\geq 0$ for every $v_1,v_2\in X$.
		\item $A$ is maximal monotone if $A$ is monotone and if $v_1\in X$ and $b\in X^*$ such that $\langle b-A(v_2),v_1-v_2\rangle \geq 0$ for every $v_2\in X$, then $A(v_1)=b$.
		\item $A$ is pseudomonotone if for every sequence $(v_k)_k$ in $X$ converging to $v$ in $X$ and $\displaystyle\limsup_{k\rightarrow\infty}\langle A(v_k),v_k-v\rangle\leq 0$, it follows that $(A(v_k))_k$ converges to $A(v)$ weakly in $X^{*}$ and $\displaystyle\lim_{k\rightarrow\infty}\langle A(v_k),v_k-v\rangle=0$.
		\item $A$ is coercive if $\displaystyle\lim_{||v||_{X}\rightarrow\infty}\frac{\langle A(v),v\rangle}{||v||_X}=+\infty$.
	\end{itemize}
Further, let a linear operator $L$ be maximal monotone and densely defined map from the domain $D(L)\subset X$ to $X^*$. Then we say that a bounded, demicontinuous map $A:X\rightarrow X^*$ is pseudomonotone with respect to $D(L)$ if for any sequence $(v_k)_k$ in $D(L)$ with $v_k\rightarrow v$ weakly in $X$, $Lv_k\rightarrow Lv$ weakly in $X^*$ and $\displaystyle\limsup_{k\rightarrow\infty}\langle A(v_k),v_k-v\rangle\leq 0$, we have $\displaystyle\lim_{n\rightarrow\infty}\langle A(v_k),v_k-v\rangle=0$ and $A(v_k)\rightarrow A(v)$ weakly in $X^*$.
\end{definition}

Following \cite{lsimon},  we define the differential operator $ L: V\rightarrow W^*$   as
	\begin{equation}\label{diffop}
		Lv:=\partial_tv
	\end{equation}
and by $D(L)$, we denote the domain of the map $L$. We have $D(L)=V$.
The operator $L$ is maximal monotone, closed and densely defined linear operator \cite{lsimon}.\\
 Now, we recall the results which are used to show the existence of the weak solution of problem \eqref{eqn1.1}.
Consider the following problem,
\begin{equation}\label{eqngen}
	\partial_tu-\sum_{i=1}^{d}\frac{\partial}{\partial x_i}(a_i(\boldsymbol{x},t,u(\boldsymbol{x},t),\nabla u(\boldsymbol{x},t);u))+a_0(\boldsymbol{x},t,u(\boldsymbol{x},t),\nabla u(\boldsymbol{x},t);u)=f,
\end{equation}
where $a_i$ are real valued functions defined on $\Omega\times(0,T)\times\mathbb R^{d+1}\times W$. Problem \eqref{eqn1.1} can be considered as a special case of \eqref{eqngen}. In \cite{lsimon}, the following theorems (Theorems \ref{thm2.2} and \ref{thm2.1}) were used to prove the existence of the weak solution of \eqref{eqngen}. 

\begin{theorem}\label{thm2.2}\cite[Theorem 5.1]{lsimon}
Let $A:W\rightarrow W^*$ be the operator defined as 
\begin{equation*}
\begin{split}
    &\langle A(u),v\rangle=\\
    &\int_{0}^{T}\int_{\Omega}\left[\sum_{i=1}^{d}a_i(\boldsymbol{x},t,u(\boldsymbol{x},t),\nabla u(\boldsymbol{x},t);u)\frac{\partial v}{\partial x_i}+a_0(\boldsymbol{x},t,u(\boldsymbol{x},t),\nabla u(\boldsymbol{x},t);u)v\right]\ d\Omega\ dt
\end{split}
\end{equation*} 
$\forall u,v\in W.$\\
Assume that the functions $a_i$ satisfy the following conditions:
\begin{enumerate}[(C1)]
	\item The functions $a_i:\Omega\times(0,T)\times\mathbb R^{d+1}\times W\rightarrow\mathbb R$ satisfy the Carath\'eodory conditions for fixed $u\in W$, i.e., they are measurable in $(\x, t)$ for each fixed $(\eta_0,\boldsymbol{\eta})\in\mathbb{R}^{d+1}$ and continuous in $(\eta_0,\boldsymbol{\eta})$ for a.e. fixed $(\x,t)\in \Omega\times (0,T)$ and they have the Volterra property, i.e., $a_i(\x,t,\boldsymbol{\eta}_0,\boldsymbol{\eta};u)$ depends only on the restriction of $u$ to $[0,t]$ for $i=0,\ldots,d$.
	\item There exists bounded operators $h_1:W\rightarrow \mathbb R^+$and $l_1:W\rightarrow Z$ such that
	\begin{equation*}
		\vert a_i(\boldsymbol{x},t,\eta_0,\boldsymbol{\eta};u)\vert\leq h_1(u)[\vert\eta_0\vert+\vert\boldsymbol{\eta}\vert]+[l_1(u)](\boldsymbol{x},t),\ i=0,1,\cdots, d,
	\end{equation*}
	for a.e. $(\boldsymbol{x},t)\in\Omega\times (0,T)$, each $(\eta_0,\boldsymbol{\eta})\in\mathbb R^{d+1}$, with $\boldsymbol{\eta}\in\mathbb{R}^d$ and $u\in W$.
	\item $\sum_{i=1}^{d}[a_i(\boldsymbol{x},t,\eta_0,\boldsymbol{\eta};u)-a_i(\boldsymbol{x},t,\eta_0,\boldsymbol{\eta}^*;u)](\eta_i-\eta_i^*)>0$ if $\boldsymbol{\eta}\neq \boldsymbol{\eta}^*$.
	\item There exist bounded operators $h_2:W\rightarrow \mathbb R^+$, $l_2:W\rightarrow L^1(\Omega\times(0,T))$ such that
	\begin{equation*}
		\sum_{i=0}^{d}a_i(\boldsymbol{x},t,\eta_0,\boldsymbol{\eta};u)\eta_i\geq [h_2(u)][\vert\eta_0\vert^2+\vert\boldsymbol{\eta}\vert^2]-[l_2(u)](\boldsymbol{x},t)
	\end{equation*}
	for a.e. $(\boldsymbol{x},t)\in\Omega\times(0,T)$, all $(\eta_0,\boldsymbol{\eta})\in\mathbb R^{d+1}$,  with $\boldsymbol{\eta}\in\mathbb{R}^d$,  $u\in W$ and
	\begin{equation*}
		\displaystyle \lim_{\vert\vert u\vert\vert_W\rightarrow\infty}\left[h_2(u)\vert\vert u\vert\vert_W-\frac{\vert\vert l_2(u)\vert\vert_{L^1(\Omega\times(0,T))}}{\vert\vert u\vert\vert_W}\right]=+\infty.
	\end{equation*}
	\item There exists $\theta>0$ such that if $(u_n)$ converges to $ u$ weakly in $W$, strongly in $L^2(0,T;H^{1-\theta}(\Omega))$ then for $i=0,1,\cdots,d$
	\begin{equation*}
		a_i(\boldsymbol{x},t,u_n(\boldsymbol{x},t),\nabla u_n(\boldsymbol{x},t);u_n)-a_i(\boldsymbol{x},t,u_n(\boldsymbol{x},t),\nabla u_n(\boldsymbol{x},t);u)\rightarrow 0
	\end{equation*}
	in $Z$.
\end{enumerate}
Then the operator $A:W\rightarrow W^*$ is  bounded, demicontinuous, pseudomonotone with respect to $V$, coercive and of Volterra type.
\end{theorem}

\begin{theorem}\label{thm2.1}\cite[Theorem 2.4]{lsimon}
	If $X$ is a reflexive Banach space, $L:X\supseteq D(L)\rightarrow X^*$ is a closed, maximal monotone, densely defined  linear operator and $A:X\rightarrow X^*$ is  bounded, coercive, demicontinuous and pseudomonotone with respect to $D(L)$, then $(L+A)(D(L))=X^*$.
\end{theorem}

	\section{Space-time weak formulation: existence-uniqueness results}\label{sec2}
	Multiplying the parabolic equation \eqref{eqn1.1} with a test function $v\in W$, integrating  over $\Omega\times(0,T)$ and integrating by parts in the nonlinear elliptic term, the space-time weak formulation of problem \eqref{eqn1.1} reads: Find $u\in V$ such that
	\begin{equation}\label{weak1}
		\int_{0}^{T}\int_{\Omega} \partial_tu v\ d\Omega\ dt+\int_{0}^{T}\int_{\Omega} a(l(u))\nabla u\cdot\nabla v\ d\Omega\ dt=\int_{0}^{T}\int_{\Omega}fv\ d\Omega\ dt\quad \forall v\in W.
	\end{equation}
Using the distribution theory, one can show that the above formulation is equivalent to the following weak formulation: Find $u\in V$ such that
\begin{equation}\label{linevf}
    \langle \partial_tu,v\rangle+a(l(u))(\nabla u,\nabla v)=\langle f,v \rangle\quad \forall v\in H_0^1(\Omega)\quad \text{for a.e. } t\in(0,T).
\end{equation}
Note that the existence-uniqueness of weak solution corresponding to weak formulation \eqref{linevf} is already established in \cite{Chipot1999} using Faedo-Galerkin method. However, here we give a different proof for the existence of weak solution corresponding to space-time formulation \eqref{weak1} using the theory of monotone type of operators.
 In our further analysis, we consider the following hypotheses:
	\begin{enumerate}[(H1)]
		\item There exists two positive real numbers $m$ and $M$ such that,
		\begin{equation*}
			0< m\leq a(s)\leq M<\infty\ \forall s\in\mathbb R.
		\end{equation*}
		\item  The function $a$ is Lipschitz continuous with Lipschitz constant $L_M$, i.e., for any $s_1,s_2\in\mathbb{R}$, the following holds
        \begin{equation*}
            \vert a(s_1)-a(s_2)\vert\leq L_M\vert s_1-s_2\vert.
        \end{equation*}
	\end{enumerate}

Now, in the following theorem, we prove the existence and uniqueness of the weak solution of the problem \eqref{eqn1.1}.
    \begin{theorem}\label{thm3.1}
	 Let $f\in W^*$ and assume that  (H1)-(H2) hold. Then there exists a unique $u\in V$ which satisfies problem \eqref{weak1}.
	\end{theorem}
\begin{proof}
	In order to prove the existence of a solution of \eqref{weak1}, we define a nonlinear operator $A:W\rightarrow W^*$ by
	\begin{equation}\label{A}
		\langle A(u),v\rangle:=\int_{0}^{T}\int_{\Omega} a(l(u))\nabla u\cdot\nabla v\ d\Omega\ dt\quad \forall v\in W.
	\end{equation}
	Then $A$ is coercive, bounded, demicontinuous and pseudomonotone with respect to $V$. The proof follows from Theorem \ref{thm2.2} by taking $a_0(\boldsymbol{x},t,\eta_0,\boldsymbol{\eta};u)=0$ and $a_i(\boldsymbol{x},t,\eta_0,\boldsymbol{\eta};u)=a(l(u))\eta_i$ for $i=1,\cdots,d$.
	Since the differential operator $L$ defined in \eqref{diffop} is densely defined, closed and maximal monotone operator, from Theorem \ref{thm2.1} it follows that there exists  $u\in V$ which satisfies \eqref{weak1} .\\
To show the uniqueness of the weak solution, we follow the idea given in \cite{ROBALO20145609,Glazyrina2012}.
	\\ Suppose $u_1,u_2$ are two solutions of problem \eqref{weak1}. From \eqref{weak1},   we have
	\begin{equation}\label{weaksol1}
		\langle \partial_tu_1(t),v(t)\rangle+ (a(l(u_1))\nabla  u_1(t),\nabla  v(t))=\langle f(t),v(t)\rangle
	\end{equation}
	and 
	\begin{equation}\label{weaksol2}
		\langle \partial_tu_2(t),v(t)\rangle+ (a(l(u_2))\nabla  u_2(t),\nabla  v(t))=\langle f(t),v(t)\rangle
	\end{equation}for almost all $t\in(0,T)$.\\
 Subtracting \eqref{weaksol2} from \eqref{weaksol1} and defining
	$\psi:=u_1-u_2$, we get
	\begin{equation*}
		\begin{split}
			&\langle \partial_t\psi,\psi\rangle+a(l(u_1))(\nabla  u_1,\nabla  \psi)-a(l(u_1))(\nabla  u_2,\nabla  \psi) \\
			&= a(l(u_2))(\nabla  u_2,\nabla  \psi)-a(l(u_1))(\nabla  u_2,\nabla  \psi).
		\end{split}
	\end{equation*}
	Then,
	\begin{equation*} 
			\displaystyle\frac{1}{2}\frac{d}{dt}||\psi(t)||_{L^2(\Omega)}^2+a(l(u_1))||\nabla  \psi(t)||_{L^2(\Omega)}^2=  (a(l(u_2))-a(l(u_1)))(\nabla  u_2,\nabla  \psi). 
	\end{equation*}
	We also have $a(l(u_1))||\nabla  \psi(t)||_{L^2(\Omega)}^2\geq m||\nabla  \psi(t)||_{L^2(\Omega)}^2$ and as $a$ is a Lipschitz continuous function, we have
	\begin{equation*}
		\begin{split}
			|a(l(u_1))-a(l(u_2))| & \leq L_M|l(u_1)-l(u_2)| \\
			& \leq C||u_1(t)-u_2(t)||_{L^2(\Omega)}.
		\end{split}
	\end{equation*}
	Consequently,
	\begin{equation*}
		\displaystyle\frac{d}{dt}||\psi(t)||_{L^2(\Omega)}^2+2m||\nabla  \psi(t)||_{L^2(\Omega)}^2\leq C||\psi(t)||_{L^2(\Omega)}||\nabla  u_2(t)||_{L^2(\Omega)}||\nabla  \psi(t)||_{L^2(\Omega)}.
	\end{equation*}
	Thus,  \begin{equation*}
		\displaystyle\frac{d}{dt}||\psi(t)||_{L^2(\Omega)}^2\leq C(t)||\psi(t)||_{L^2(\Omega)}^2.
	\end{equation*}
	Then, by using the Gronwall's lemma, we get
	\begin{equation*}
		||\psi(t)||_{L^2(\Omega)}^2\leq 0.
	\end{equation*}
	Thus, $\psi(t)=0$ and $u_1=u_2$.
    
\end{proof}
\section{Space-time isogeometric discretization}\label{sec4}
	For the space-time isogeometric discretization of the weak formulation \eqref{weak1}, we consider the isogeometric spaces $V_{h}$ defined in \eqref{eq:isospace} and, following \cite{LOLI20202586}, we take $W_h=V_{h}$. Then the isogeometric approximation of \eqref{weak1} is given by: Find $u_h\in V_h$ such that
	\begin{equation}\label{eqn5.1}
		\int_{0}^{T}\int_{\Omega}\partial_t u_{h}v_h\ d\Omega\   dt+\int_{0}^{T}\int_{\Omega}^{}a(l(u_h))\nabla u_h\cdot\nabla  v_h\ d\Omega\   dt=\int_{0}^{T}\int_{\Omega}^{}fv_h\ d\Omega\   dt\quad \forall v_h\in V_h.
	\end{equation}
 	 
	Now, in the following lemma, we prove a priori estimate of the solution $u_h$ of problem \eqref{eqn5.1}.
    \begin{lemma}\label{aprioriestimate}
	Suppose that (H1)-(H2) hold and $f\in Z$. Then any solution $u_h$ of \eqref{eqn5.1} satisfies the following		
        \begin{equation*}
			\vert\vert u_h\vert\vert_{W}\leq \frac{{C_p}}{m}\vert\vert f\vert\vert_Z.
		\end{equation*}
		
	\end{lemma}

	\begin{proof} 
We take $v_h=u_h$ in \eqref{eqn5.1} and get
\begin{equation*}
 \int_{0}^{T}\int_{\Omega}\partial_t u_{h}u_h\ d\Omega\   dt+\int_{0}^{T}\int_{\Omega}^{}a(l(u_h))\nabla u_h\cdot\nabla  u_h\ d\Omega\   dt=\int_{0}^{T}\int_{\Omega}^{}fu_h\ d\Omega\ dt.   
\end{equation*}
Then
\begin{equation*}
    \frac{1}{2}\vert\vert u_h(T)\vert\vert_{L^2(\Omega)}^2+\int_{0}^{T}\int_{\Omega}a(l(u_h))\nabla u_h\cdot\nabla u_h\ d\Omega dt\leq \vert\vert f\vert\vert_Z\vert\vert u_h\vert\vert_Z.
\end{equation*}
Using the assumption on $a$ and Poincar\'e inequality, we get
\begin{equation*}
    m\vert\vert u_h\vert\vert_W^2\leq C_p\vert\vert f\vert\vert_Z\vert\vert u_h\vert\vert_W.
\end{equation*}
Finally, we have
\begin{equation*}
    \vert\vert u_h\vert\vert_{W}\leq \frac{{C_p}}{m}\vert\vert f\vert\vert_Z.
\end{equation*}
\end{proof}
Next, in the following, we prove the existence of a solution to the problem \eqref{eqn5.1}. We follow the main ideas of \cite[Section 3]{gudi2012finite}, where similar estimates have been proven for a stationary nonlocal problem and utilize a consequence of the Brouwer's fixed point theorem, which is stated below.
\begin{lemma}\cite{Kesavan}\label{brouwer}
Let $H$ be a finite dimensional Hilbert space equipped with inner product $(\cdot,\cdot)_H$ and norm $\vert\cdot\vert_H$. Let $R:H\rightarrow H$ be a continuous map satisfying the following property: there exists $\rho>0$ such that
\begin{equation*}
    (R(w),w)>0\quad\mbox{for all }w\in H\mbox{ with }\vert w\vert_{H}=\rho.
\end{equation*}
Then, there exists  $u\in H$ such that $R(u)=0$ and $\vert u\vert_H\leq \rho$ .
\end{lemma}

    \begin{theorem}\label{exisproof}
	Let (H1)-(H2) hold and $f\in Z$. Then there exists a solution $u_h\in V_h$ of \eqref{eqn5.1}.
	\end{theorem}
	\begin{proof} We consider $V_h$ as a finite dimensional subspace of the Hilbert space $W$ with inner product
	\begin{equation*}
		(u_h,v_h)_W=\int_{0}^{T}\int_{\Omega}\nabla  u_h\cdot\nabla  v_h\ d\Omega\ dt.
	\end{equation*}
    Fix $w_h\in V_h$. Define $Q_{w_h}:V_h\rightarrow \mathbb{R}$ by
	\begin{equation*}
		Q_{w_h}(v_h)=\int_{0}^{T}\int_{\Omega}\partial_tw_hv_h\ d\Omega\ dt+\int_{0}^{T}\int_{\Omega}a(l(w_h))\nabla  w_h\cdot\nabla  v_h\ d\Omega\ dt-\int_{0}^{T}\int_{\Omega}fv_h \ d\Omega\ dt.
	\end{equation*}
	Then, \begin{equation*}
		\begin{split}
			\vert Q_{w_h}(v_h)\vert&\leq \vert\vert \partial_tw_h\vert\vert_Z\vert\vert v_h\vert\vert_Z+ M\vert\vert\nabla  w_h\vert\vert_Z\vert\vert\nabla  v_h\vert\vert_Z+\vert\vert f\vert\vert_Z\vert\vert v_h\vert\vert_Z\\
			&\leq (C\vert\vert w_h\vert\vert_W+C_p\vert\vert f\vert\vert_Z)\vert\vert v_h\vert\vert_W,
		\end{split}
	\end{equation*}
    where we used that $\|\partial_t w_h\|^2_Z\leq \|\partial_t w_h\|^2_Z+\|  w_h\|^2_W\leq \| w_h\|^2_V \leq C^2 \|w_h\|^2_W$.
Thus, $Q_{w_h}$ is a continuous linear map and hence by the Riesz representation theorem, we get $R(w_h)\in V_h$ such that
\begin{equation*}
	(R(w_h),v_h)_W=Q_{w_h}(v_h)\ \forall v_h\in V_h.
\end{equation*} 
We show that the function $R:V_h\rightarrow V_h$ is continuous.\\
Fix $w_h\in V_h$ and let $\epsilon>0$. Then for any $v_h,\phi_h\in V_h$, we have
\begin{equation*}
	\begin{split}
		(R(w_h)-R(\phi_h),v_h)_W=&(R(w_h),v_h)_W-(R(\phi_h),v_h)_W\\
						=&\int_{0}^{T}\int_{\Omega}(\partial_tw_h-\partial_t\phi_h)v_h\ d\Omega\ dt\\
      &+\int_{0}^{T}\int_{\Omega}(a(l(w_h))\nabla  w_h-a(l(\phi_h))\nabla  \phi_h)\cdot\nabla  v_h\ d\Omega\   dt.
	\end{split}
\end{equation*}
Then, using the fact that $a$ is Lipschitz and the embedding $V\hookrightarrow C(0,T;L^2(\Omega))$, we get
\begin{equation*}
	\begin{split}
		&\int_{0}^{T}\int_{\Omega}(a(l(w_h))\nabla  w_h-a(l(\phi_h))\nabla  \phi_h)\cdot\nabla  v_h\ d\Omega\ dt\\
		&=\int_{0}^{T}\int_{\Omega}(a(l(w_h))-a(l(\phi_h)))\nabla  w_h\cdot\nabla  v_h\ d\Omega\ dt+\int_{0}^{T}\int_{\Omega}a(l(\phi_h))(\nabla  w_h-\nabla \phi_h)\cdot\nabla  v_h\ d\Omega\ dt\\
		&\leq C'\vert\vert \nabla  w_h-\nabla \phi_h\vert\vert_Z\vert\vert\nabla  w_h\vert\vert_Z\vert\vert \nabla  v_h\vert\vert_Z+M\vert\vert \nabla  w_h-\nabla \phi_h\vert\vert_Z\vert\vert\nabla  v_h\vert\vert_Z.
	\end{split}
\end{equation*}
Thus,
\begin{equation*}
	(R(w_h)-R(\phi_h),v_h)_W\leq (C+C'\vert\vert w_h\vert\vert_W+M) \vert\vert v_h\vert\vert_W\vert\vert w_h-\phi_h\vert\vert_W.
\end{equation*}
Replacing $v_h$ by $R(w_h)-R(\phi_h)$, we get
\begin{equation*}
\begin{split}
   & \vert\vert R(w_h)-R(\phi_h)\vert\vert^2_W\leq (C+C'\vert\vert w_h\vert\vert_W+M)	\vert\vert R(w_h)-R(\phi_h)\vert\vert_W	\vert\vert w_h-\phi_h\vert\vert_W.
\end{split}
\end{equation*}
Hence \begin{equation*}
		\vert\vert R(w_h)-R(\phi_h)\vert\vert_W\leq (C+C'\vert\vert w_h\vert\vert_W+M)		\vert\vert w_h-\phi_h\vert\vert_{W}.
\end{equation*}
Taking $\delta=\frac{\epsilon}{(C+C'\vert\vert w_h\vert\vert_W+M)}$, the continuity of $R$ is proved.\\\\
Now for $w_h\in V_h$,
\begin{equation*}
	\begin{split}
		(R(w_h),w_h)_W&=\int_{0}^{T}\int_{\Omega}\partial_tw_hw_h\ d\Omega\ dt+\int_{0}^{T}\int_{\Omega}a(l(w_h))\nabla  w_h\cdot\nabla  w_h\ d\Omega\ dt-\int_{0}^{T}\int_{\Omega}fw_h\ d\Omega\ dt\\
		&\geq(m\vert\vert w_h\vert\vert_W-C_p\vert\vert f\vert\vert_Z)\vert\vert w_h\vert\vert_W.
	\end{split}
\end{equation*}
Then  $(R(w_h),w_h)_W>0$ for $w_h\in V_h$ having the property $\vert\vert w_h\vert\vert_W=\rho$, where $\rho$ is any number such that
\begin{equation*}
	\rho>\frac{C_p\vert\vert f\vert\vert_Z}{m}.
\end{equation*}
Then the existence of the solution $u_h$ of \eqref{eqn5.1} follows from Lemma \ref{brouwer}.
\end{proof}
\begin{remark}
In Theorem \ref{exisproof}, we have proven the existence of a solution of \eqref{eqn5.1}. However, the proof of the uniqueness for the discrete solution $u_h$ is yet to be discovered.
\end{remark}
\subsection{Discrete system and linear solver}
\label{sec:4.1}

The nonlinear system associated with the discrete formulation \eqref{eqn5.1} is given by
	\begin{equation}\label{eqn5.2}
		\mathbf{A}( \mathbf{u})\mathbf{ {u}}=\mathbf{f},
	\end{equation}
	where 
	\begin{equation*}
		[\mathbf{A}( \mathbf{u})]_{j,k}=\int_{0}^{T}\int_{\Omega}( \partial_tB_{{k,\mathbf{q}}}B_{j,\mathbf{q}}+a(l(\mathbf{u}))\nabla  B_{k,\mathbf{q}}\cdot\nabla  B_{j,\mathbf{q}})\ d\Omega  dt,\quad [\mathbf{f}]_{j}=\int_{0}^{T}\int_{\Omega}fB_{j,\mathbf{q}}\  d\Omega \ dt
	\end{equation*}
	and where, with a little abuse of notation, we defined $l(\mathbf{u}):=\int_{\Omega} \sum_{i=1}^{N_{dof}}\mathbf{u}_iB_{i,\mathbf{q}} d\Omega $.
	Thanks to the tensor product structure of the space $V_h$ and the fact that $a(l( \mathbf{u}))$ does not depend on the space variables, the matrix $\mathbf{A}( \mathbf{u})$ can be written as the sum of Kronecker products of matrices 
	\begin{equation*}
		\mathbf{A}( \mathbf{u})=\mathbf{W}_t\otimes\mathbf{M}_s+\mathbf{M}_t( \mathbf{u})\otimes\mathbf{K}_s, 
	\end{equation*}
	where \begin{equation*}
		[\mathbf{W}_t]_{j,k}=\int_{0}^{T}b'_{k,q_t}(t)b_{j,q_t}(t)\ dt,\quad j,k=1,\cdots,N_t,
	\end{equation*}
	\begin{equation*}
		[\mathbf{M}_s]_{j,k}=\int_{\Omega}B_{j,\mathbf{q}_s}(\boldsymbol{x})B_{k,\mathbf{q}_s}(\boldsymbol{x})\  d\Omega ,\quad j,k=1,\cdots,N_s,
	\end{equation*}
	\begin{equation*}
		[\mathbf{M}_t(\mathbf{u})]_{j,k}=\int_{0}^{T}a(l( \mathbf{u}))b_{k,q_t}(t)b_{j,q_t}(t)\ dt,\quad j,k=1,\cdots,N_t,
	\end{equation*}
	and
	\begin{equation*}
		[\mathbf{K}_s]_{j,k}=\int_{\Omega}\nabla  B_{j,\mathbf{q}_s}(\boldsymbol{x})\cdot\nabla  B_{k,\mathbf{q}_s}(\boldsymbol{x})\  d\Omega ,\quad j,k=1,\cdots,N_s.
	\end{equation*}
	We solve the nonlinear system \eqref{eqn5.2} by using Picard's method, which reads as follows: For given initial guess ${\mathbf{u}}^0$, solve
	\begin{equation}\label{eqn5.3}
		\mathbf{A}( \mathbf{u}^{(n-1)})\mathbf{ {u}}^{(n)}=\mathbf{f},
	\end{equation} for $n=1,2,\cdots$ until \begin{equation}\label{eq:stopPic}\vert\vert \mathbf{u}^{(n+1)}-\mathbf{u}^{(n)}\vert\vert_{\infty}\leq\epsilon\end{equation} for a given tolerance $\epsilon>0$.

  We now propose a solving strategy for the linear system \eqref{eqn5.3}. For the sake of simplicity, we drop the dependence of $\mathbf{A}(\mathbf{u}^{(n-1)})$ and $\mathbf{M}_t(\mathbf{u}^{(n-1)})$ from $\mathbf{u}^{(n-1)}$, and we write, instead, $\mathbf{A}$ and $\mathbf{M}_t$. Following \cite[Section 4.2.2]{LOLI20202586}, we build a stable factorization of the time pencil $(\mathbf{W}_t,\mathbf{M}_t)$ for $n\geq0$: we find $\mathbf{U}_t,\mathbf{\Delta}_t\in\mathbb{C}^{N_t\times N_t}$, depending on $\mathbf{u}^{(n-1)}$, such that
  \begin{equation}
\label{eq:factorization}\mathbf{U}_t^*\mathbf{M}_t \mathbf{U}_t=\mathbf{1}_{N_t} \quad \text{and} \quad \mathbf{U}_t^*\mathbf{W}_t\mathbf{U}_t=\mathbf{\Delta}_t, \end{equation}
  where $\mathbf{1}_{m}\in\mathbb{R}^{m\times m}$ denotes the identity matrix and $\mathbf{\Delta}_t$ has an arrowhead structure (only the diagonal, that last row and the last column could be non-zero). Then, the linear system matrix $\mathbf{A}$ factorizes as
  $$
  \mathbf{A}=(\mathbf{U}_t^{-*}\otimes \mathbf{1}_{N_s})(\mathbf{\Delta}_t\otimes\mathbf{M}_s+\mathbf{1}_{N_t}\otimes\mathbf{K}_s)(\mathbf{U}_t^{-1}\otimes \mathbf{1}_{N_s}).
  $$ 
Differently than in \cite{LOLI20202586}, we do solve  iteratively the linear system  \eqref{eqn5.3}, but we first exploit the decomposition of the system matrix.
Thus, to compute the solution of \eqref{eqn5.3}, we use  Algorithm \ref{al:solving}.
 \begin{algorithm}
\caption{Solving Algorithm}\label{al:solving}
\hspace*{\algorithmicindent} \textbf{Input}:  Linear system matrix $\mathbf{A}(\mathbf{u}^{(n-1)})$ and right-hand side $\mathbf{f}$\\
 \hspace*{\algorithmicindent} \textbf{Output}: $\mathbf{u}^{(n)}$ solution of the system  \eqref{eqn5.3} 
 \begin{algorithmic}[1] 
\State Compute the factorization \eqref{eq:factorization};
    \State Compute $\widetilde{\mathbf{f}}:=(\mathbf{U}_t^{*}\otimes \mathbf{1}_{N_s})\mathbf{f}$;
    \State Solve $(\mathbf{\Delta}_t\otimes\mathbf{M}_s+\mathbf{1}_{N_t}\otimes\mathbf{K}_s)\widetilde{\mathbf{u}}^{(n)} =  \widetilde{\mathbf{f}};$
    \State Compute $\mathbf{u}^{(n)}:=(\mathbf{U}_t\otimes \mathbf{1}_{N_s})\widetilde{\mathbf{u}}^{(n)}.$

\end{algorithmic}
\end{algorithm}


The computational cost of Step 1 is $O(N_t^3)$, while Step 2 and Step 4 yield a computational cost of $4N_{dof}N_t$ FLOPs. 
We solve the linear system in Step 3 iteratively and we use as a preconditioner the matrix
\begin{equation}
    \label{eq:precond}
\widehat{\mathbf{P}}:=\mathbf{\Delta}_t\otimes\widehat{\mathbf{M}}_s+\mathbf{1}_{N_t}\otimes\widehat{\mathbf{K}}_s
\end{equation}
where $[\widehat{\mathbf{M}}_s]_{j,k}:=\int_{\widehat{\Omega}} \widehat{B}_{k,{q}_s}\widehat{B}_{j,q_s}\ d\widehat{{\Omega}}$ and $[\widehat{\mathbf{K}}_s]_{j,k}:=\int_{\widehat{\Omega}} \nabla\widehat{B}_{k,q_s}\cdot\nabla\widehat{B}_{j,q_s}\  d\widehat{{\Omega}}$ for $j,k=1,\dots N_s$ are the mass and stiffness matrices discretized in the parametric domain.  Thanks to the tensor product structure of the basis functions, we have that
$$
\widehat{\mathbf{M}}_s=\bigotimes_{i=1}^d\widehat{\mathbf{M}}_i \quad \text{and} \quad \widehat{\mathbf{K}}_s=\sum_{i=1}^d\widehat{\mathbf{M}}_d\otimes\dots\otimes \widehat{\mathbf{M}}_{i+1}\otimes\widehat{\mathbf{K}}_{i}\otimes\widehat{\mathbf{M}}_{i-1}\otimes\dots\otimes\widehat{\mathbf{M}}_{1},
$$
where $\widehat{\mathbf{M}}_{i}$ and  $\widehat{\mathbf{K}}_{i}$ are the one-dimensional mass and stiffness matrices, respectively. As the matrices $\widehat{\mathbf{K}}_i$ and $\widehat{\mathbf{M}}_i$ are positive definite, they admit a generalized eigendecomposition: We can find $\mathbf{U}_i$ matrix of generalized eigenvectors, and $\mathbf{\Lambda}_i$ matrix of generalized eigenvalues, such that
\begin{equation}
    \label{eq:space_dec}
    \mathbf{U}_i^T\widehat{\mathbf{M}}_i\mathbf{U}_i=\mathbf{1}_{n_{s,i}} \quad \text{and} \quad  \mathbf{U}_i^T\widehat{\mathbf{K}}_i\mathbf{U}_i= \mathbf{\Lambda}_i.
\end{equation}
Factorizing the common terms, the preconditioner \eqref{eq:precond} can be written as
\begin{equation}
    \label{eq:final_precond}
    \widehat{\mathbf{P}}= (\mathbf{1}_{N_t}\otimes \mathbf{U}^{-T}_d\otimes\dots\otimes \mathbf{U}_1^{-T})(\mathbf{\Delta}_t\otimes\mathbf{1}_{N_s}+\mathbf{1}_{N_t}\otimes \mathbf{\Lambda}) (\mathbf{1}_{N_t}\otimes \mathbf{U}^{-1}_d\otimes\dots \otimes \mathbf{U}^{-1}_1),
\end{equation}
where $\mathbf{\Lambda}:=\sum_{i=1}^d\bigotimes \mathbf{1}_{n_{s,d}}\dots\otimes \mathbf{1}_{n_{s,i+1}}\otimes\mathbf{\Lambda}_i\otimes\mathbf{1}_{n_{s,i-1}}\otimes\dots\otimes \mathbf{1}_{n_{s,1}}$.

Finally, the application of the preconditioner \eqref{eq:final_precond} is given by the solution of the following equation: Given $\mathbf{r}$, find $\mathbf{s}$ such that 
\begin{equation}
\label{eq:prec_appl}
\mathbf{\widehat{P}}\mathbf{s}=\mathbf{r}
\end{equation}
and it is performed by Algorithm \ref{al:prec}. 
 \begin{algorithm}[H]
\caption{Preconditioner setup and application}\label{al:prec}
\hspace*{\algorithmicindent} \textbf{Input}:  right-hand side $\mathbf{r}$\\
 \hspace*{\algorithmicindent} \textbf{Output}: $\mathbf{s}$ solution of the system  \eqref{eq:prec_appl}
 \begin{algorithmic}[1] 
  \State Compute the factorization \eqref{eq:space_dec};
  \State  Compute $\widetilde{\mathbf{r}}:= (\mathbf{1}_{N_t}\otimes \mathbf{U}^{T}_d\otimes\dots \otimes \mathbf{U}^{T}_1) \mathbf{r}$;
    \State  Compute $\widetilde{\mathbf{s}}:= (\mathbf{\Delta}_t\otimes\mathbf{1}_{N_s}+\mathbf{1}_{N_t}\otimes \mathbf{\Lambda})^{-1} \widetilde{\mathbf{r}}$;
   \State  Compute $ {\mathbf{s}}:= (\mathbf{1}_{N_t}\otimes \mathbf{U}_d\otimes\dots \otimes \mathbf{U}_1)  \mathbf{\widetilde{s}}$.
\end{algorithmic}
\end{algorithm} 
 We briefly analyze the computational cost of the preconditioner. We suppose for simplicity that all the spatial directions have the same number of degrees of freedom, i.e. $n_{s,l} = n_s$ for $l=1,\dots,d$.  The setup of the preconditioner is performed by step 1 in Algorithm \ref{al:prec} and has a computational cost of $O(dn_s^3)$ FLOPs, which is optimal for $d=2$ and negligible for $d\geq 3$. Step 2 and step 4 require a computational cost of $4N_{dof}dn_s$, while step 2 can be easily computed exploiting an LU decomposition, yielding a computational cost of $O(N_{dof})$ (see \cite{LOLI20202586}  for details).  
 
 In our computational tests, we use a variant of the preconditioner \eqref{eq:final_precond}, that includes some information on the geometrical mapping $\widehat{\mathbf{F}}$ and the coefficients in the spatial matrices $\widehat{\mathbf{K}}_s$ and $\widehat{\mathbf{M}}_s$. 
 The idea is to make a separable approximation of the coefficients of the geometry and of the problem and incorporate the result in the stiffness and mass matrices build in the univariate spline spaces. Then, a diagonal scaling is also added. We
 use techniques similar to the ones described in \cite[Section 4.4]{LOLI20202586}. We name such modified preconditioner $\widehat{\mathbf{P}}_G$. The computational cost of the setup and application of $\widehat{\mathbf{P}}_G$ is the same as the ones of $\widehat{\mathbf{P}}$. We remark that the setup step has to be done only once in each non-linear iteration.

\section{Convergence analysis}\label{sec5}
In this section, we derive an a priori error estimate for the numerical scheme \eqref{eqn5.1}. For deriving the error estimate, we will use a result proven in  \cite[Theorem 4.1]{chaudhary2017finite}, that is the exact solution $u$ satisfies
\begin{equation}\label{bound}
    \vert\vert u\vert\vert_{L^{\infty}(0,T;H_0^1(\Omega))}<R_1,
\end{equation}
where $R_1$ is a positive constant depending on the given data $f$.
\par We begin with proving the following Cea's type lemma.

\begin{lemma}\label{cea}
	Suppose that (H1)-(H2) hold and the Lipschitz constant $L_M$ in (H2) is such that $(m-L_MC_pR_1\vert\Omega\vert^{\frac{1}{2}})>0$. Then the solution $u_h$ of \eqref{eqn5.1} satisfies
	\begin{equation}\label{ceain}
		\vert\vert u-u_h\vert\vert_W\leq C\vert\vert u-v_h\vert\vert_V\quad \forall v_h\in V_h,
	\end{equation}	
	where $C$ is a constant independent of meshsize $h$.
\end{lemma}
\begin{proof} Let $v_h\in V_h$. Then we have
\begin{equation*}
	\vert\vert u-u_h\vert\vert_W\leq \vert\vert u-v_h\vert\vert_W+	\vert\vert v_h-u_h\vert\vert_W.
\end{equation*}
Let $\phi_h=v_h-u_h$.\\
Then, we have that 

\begin{equation*}
    \begin{split}
      &\int_{0}^{T}\int_{\Omega}\partial_t \phi_h\phi_h\ d\Omega\ dt+\int_{0}^{T}\int_{\Omega}a(l(u_h))\nabla \phi_h\cdot\nabla \phi_h\ d\Omega\ dt\\
    &=\int_{0}^{T}\int_{\Omega}(\partial_tv_{h}-\partial_tu_{h})\phi_h\ d\Omega\ dt+\int_{0}^{T}\int_{\Omega}a(l(u_h))\nabla (v_h-u_h)\cdot\nabla \phi_h\ d\Omega\ dt\\
    &=\int_{0}^{T}\int_{\Omega}(\partial_tv_{h}-\partial_tu)\phi_h\ d\Omega\ dt-\int_{0}^{T}\int_{\Omega}a(l(u))\nabla u\cdot\nabla \phi_h\ d\Omega\ dt\\
    &+\int_{0}^{T}\int_{\Omega}a(l(u_h))\nabla v_h\cdot\nabla \phi_h\ d\Omega\ dt
    \end{split}
    \end{equation*}
Adding and subtracting $\int_{0}^{T}\int_{\Omega}a(l(u_h))\nabla u\cdot\nabla \phi_h\ d\Omega\ dt$ in above, we get
    \begin{equation*}
    \begin{split}
    &\int_{0}^{T}\int_{\Omega}\partial_t \phi_h\phi_h\ d\Omega\ dt+\int_{0}^{T}\int_{\Omega}a(l(u_h))\nabla \phi_h\cdot\nabla \phi_h\ d\Omega\ dt\\
    &=\int_{0}^{T}\int_{\Omega}(\partial_tv_{h}-\partial_tu)\phi_h\ d\Omega\ dt+\int_{0}^{T}\int_{\Omega}a(l(u_h))(\nabla  v_h-\nabla  u)\cdot\nabla  \phi_h\ d\Omega\ dt\\
     &+\int_{0}^{T}\int_{\Omega}(a(l(u_h))-a(l(u)))\nabla  u\cdot\nabla \phi_h\ d\Omega\ dt\\
     &= I_1+I_2+I_3.
    \end{split}
\end{equation*}
We now estimate each term $I_1,\ I_2,\ I_3$.\\
We have
\begin{equation*}
	\begin{split}
		I_1=\int_{0}^{T}\int_{\Omega}(\partial_tv_{h}-\partial_tu)\phi_h\ d\Omega\   dt \leq \vert\vert v_h-u\vert\vert_V\vert\vert\phi_h\vert\vert_W.
	\end{split}
\end{equation*}
For $I_2$, we have
\begin{equation*}
	\begin{split}
		I_2= \int_{0}^{T}\int_{\Omega}a(l(u_h))(\nabla  v_h-\nabla  u)\cdot\nabla \phi_h\ d\Omega\ dt \leq   M\vert\vert v_h-u\vert\vert_V\vert\vert\phi_h\vert\vert_W.
	\end{split}
\end{equation*}
And for the term $I_3$, we work as follows,
\begin{equation*}
	\begin{split}
		I_3&=\int_{0}^{T}\int_{\Omega}(a(l(u_h))-a(l(u)))\nabla  u\cdot\nabla  \phi_h\ d\Omega\ dt\\
		&=\int_{0}^{T}\int_{\Omega}(a(l(u_h))-a(l(v_h)))\nabla  u\cdot\nabla  \phi_h\ d\Omega\ dt\\
		&+\int_{0}^{T}\int_{\Omega}(a(l(v_h))-a(l(u)))\nabla  u\cdot\nabla  \phi_h\ d\Omega\ dt\\
		&=I_{31}+I_{32}.
	\end{split}
\end{equation*}
Now,  using the Lipschitz continuity of $a$, the a priori estimate \eqref{bound} of $u$ and Poincar\'e inequality, we have
\begin{equation*}
	\begin{split}
		I_{31}=&\int_{0}^{T}\int_{\Omega}(a(l(u_h))-a(l(v_h)))\nabla  u\cdot\nabla  \phi_h\ d\Omega\ dt\\
		&\leq \int_{0}^{T}|(a(l(u_h))-a(l(v_h)))|\vert\vert\nabla  u(t)\vert\vert_{L^2(\Omega)}\vert\vert\nabla \phi_h(t)\vert\vert_{L^2(\Omega)}\ dt\\
		&\leq L_MC_pR_1\vert\Omega\vert^{\frac{1}{2}}\vert\vert\phi_h\vert\vert_W^2.
	\end{split}
\end{equation*}
Similarly,
\begin{equation*}
	I_{32}\leq L_MC_pR_1\vert\Omega\vert^{\frac{1}{2}}\vert\vert\phi_h\vert\vert_W\vert\vert v_h-u\vert\vert_W.
\end{equation*}
Thus, we have
\begin{equation*}
	\begin{split}
		\frac{1}{2}\vert\vert\phi_h(T)\vert\vert_{L^2(\Omega)}^2+&\int_{0}^{T}\int_{\Omega}a(l(u))\nabla \phi_h\cdot\nabla \phi_h\ d\Omega\ dt\\
		\leq & \vert\vert v_h-u\vert\vert_V\vert\vert\phi_h\vert\vert_W+M\vert\vert v_h-u\vert\vert_V\vert\vert\phi_h\vert\vert_W\\
		&+L_MC_pR_1\vert\Omega\vert^{\frac{1}{2}}\vert\vert\phi_h\vert\vert^2_W+L_MC_pR_1\vert\Omega\vert^{\frac{1}{2}}\vert\vert v_h-u\vert\vert_V\vert\vert\phi_h\vert\vert_W.
	\end{split}
\end{equation*}
Hence,
\begin{equation*}
    m\vert\vert\phi_h\vert\vert^2_W
		\leq (1+M+L_MC_pR_1\vert\Omega\vert^{\frac{1}{2}})\vert\vert v_h-u\vert\vert_V\vert\vert\phi_h\vert\vert_W+L_MC_pR_1\vert\Omega\vert^{\frac{1}{2}}\vert\vert\phi_h\vert\vert^2_W.
\end{equation*}
Then
\begin{equation*}
	(m-L_MC_pR_1\vert\Omega\vert^{\frac{1}{2}})\vert\vert\phi_h\vert\vert_W\leq C\vert\vert v_h-u\vert\vert_V,
\end{equation*}
where $C$ is a constant independent of the meshsize $h$.
And hence
\begin{equation*}
	\vert\vert\phi_h\vert\vert_W\leq C\vert\vert v_h-u\vert\vert_V.
\end{equation*}
Combining all the estimates, we get
\begin{equation*}
	\begin{split}
		\vert\vert u-u_h\vert\vert_W \leq \vert\vert u-v_h\vert\vert_W+\vert\vert\phi_h\vert\vert_W \leq C\vert\vert u- v_h\vert\vert_V,
	\end{split}
\end{equation*}
which proves the required inequality given in \eqref{ceain}.
\end{proof}
In the following theorem, we prove the a priori error estimate in $\vert\vert \cdot\vert\vert_W$ norm.

 \begin{theorem}\label{thm6.1}
Let $r$ be an integer such that $1< r\leq \min\{q_s,q_t\}+1$. If $u\in V\cap H^r(\Omega\times(0,T))$  is the solution of \eqref{weak1} and $u_h\in V_h$ the solution of \eqref{eqn5.1}, then under the assumptions of Lemma \ref{cea}, we have 
     \begin{equation}\label{eqn6.2}
		||u-u_h||_W\leq C(h_s^{r-1}+h_t^{r-1})||u||_{H^r(\Omega\times(0,T))},
	\end{equation}
 where the constant $C$ does not depend on meshsize $h$.
\end{theorem}

\begin{proof}
 Combining the Lemma \ref{cea} and the approximation estimates given in \eqref{eqprojerr}, we can prove the inequality given in \eqref{eqn6.2}. We outline the main steps here.\\
We take $v_h=\mathcal{P}_h u$ in \eqref{ceain}. So, we get
\begin{equation*}
    \vert\vert u-u_h\vert\vert_W
		\leq C\vert\vert u- \mathcal{P}_h u\vert\vert_V.
\end{equation*}
Then using the estimates of \eqref{eqprojerr}, we get the required estimates of \eqref{eqn6.2}.   
\end{proof}

\section{Numerical experiments}\label{sec7}
 \par In this section, first we present some numerical experiments to confirm the convergence estimate of Section \ref{sec5}. Then, we provide a comparison between Crank-Nicolson  time-stepping method and the space-time method. We also show some numerical results to analyze the performance of the preconditioner. All tests are performed with Matlab R2024a and GeoPDEs toolbox \cite{geopde}, on a Intel Xeon processor running at 3.90 GHz, with 256 GB RAM.
 \par In each example, we consider the final time $T=1$ and fix the tolerance $\epsilon=10^{-10}$ for stopping the Picard's iterations in \eqref{eq:stopPic}. We use the preconditioned GMRES method to solve the linear system \eqref{eqn5.3} with tolerance $10^{-12}$, the null vector as an initial guess and the preconditioner $\widehat{\mathbf{P}}_G$ described in Section \ref{sec:4.1}.  We set the meshsizes as $h_s=h_t=:h$ and the degree of the basis functions as $q_s=q_t=:q$.
   \begin{figure}
 \centering
\begin{subfigure}{0.4\textwidth}
    \includegraphics[height=7cm,width=7cm]{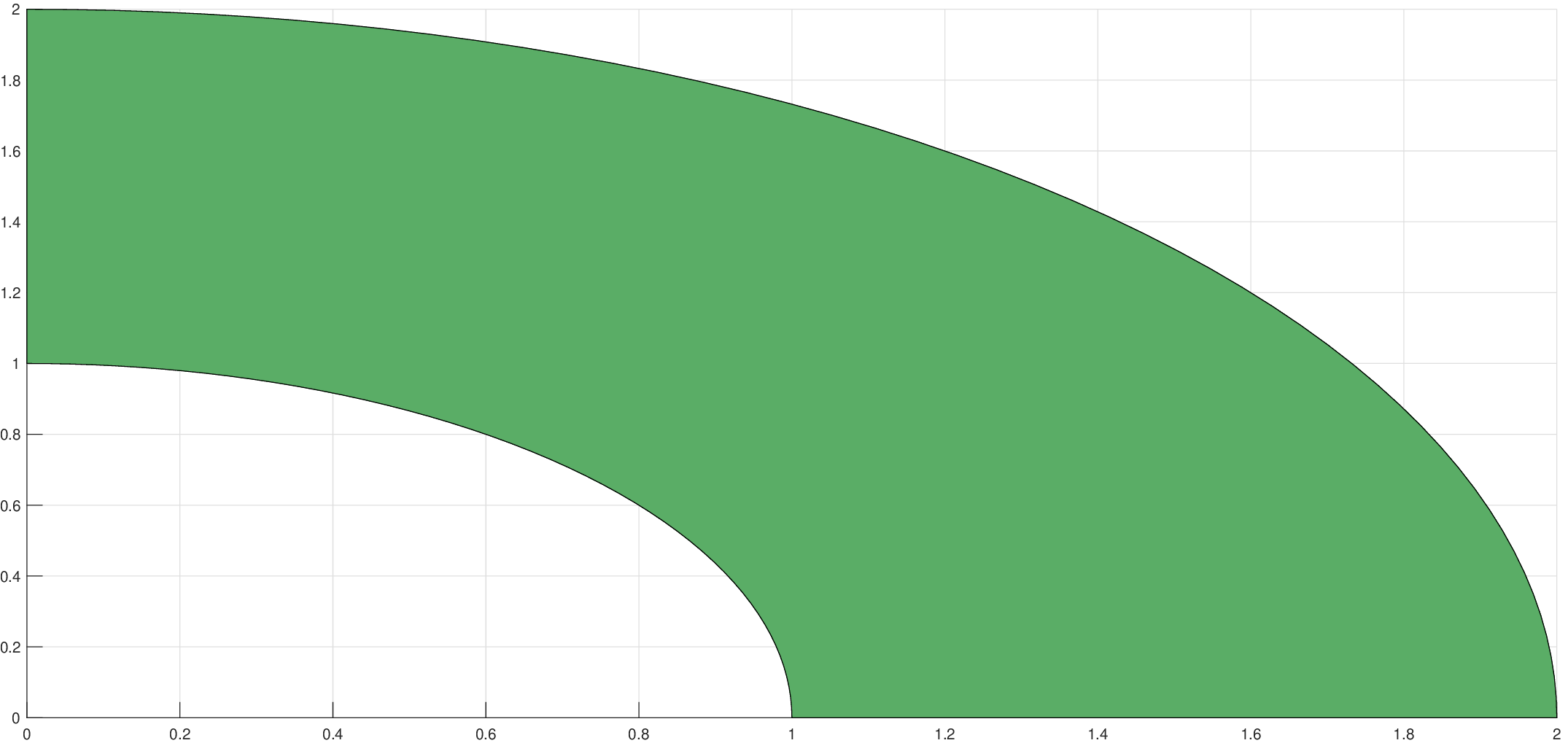}
     \caption{Quarter of annulus domain.}
    \label{fig:firstring}
\end{subfigure}
\hspace{1.5cm}
\begin{subfigure}{0.4\textwidth}
    \includegraphics[height=7cm,width=7cm]{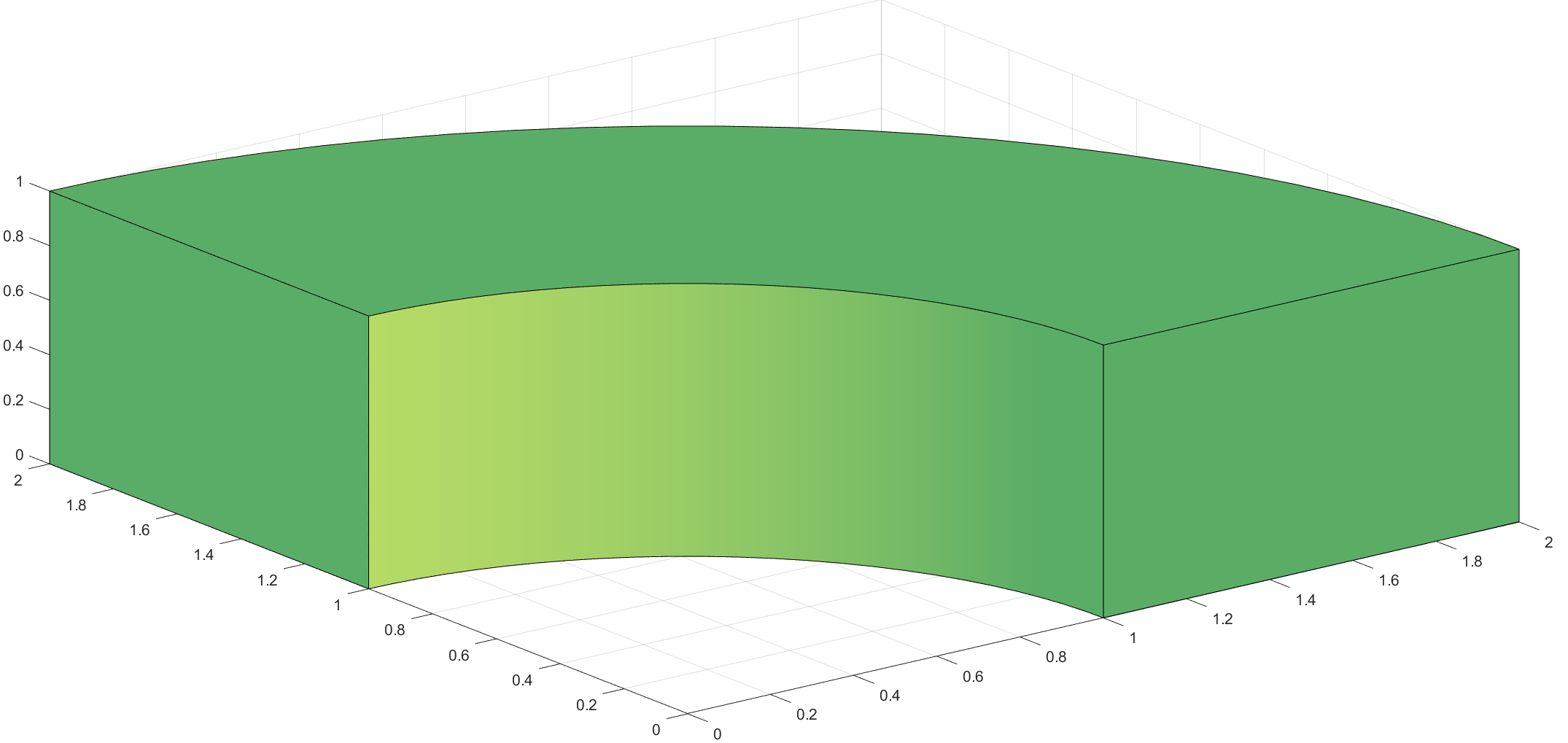}
     \caption{Thick ring shaped domain.}
    \label{fig:thick_ring}
\end{subfigure}
\hspace{1.5cm}
\begin{subfigure}{0.4\textwidth}
    \includegraphics[height=7cm,width=7cm]{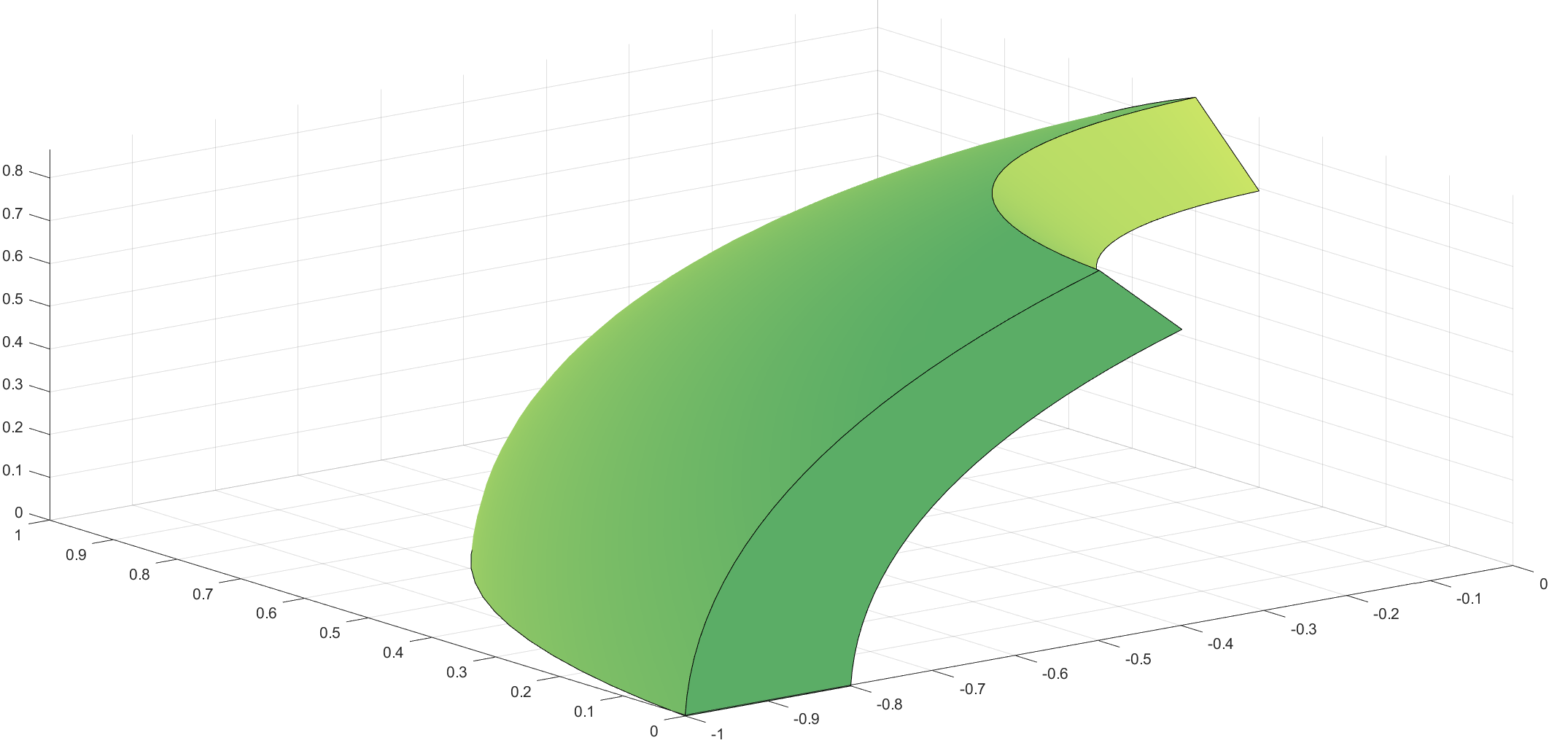}

    \caption{Igloo shaped domain.}
    \label{fig:second}
\end{subfigure}
\caption{Computational domains.}
\label{fig4}
\end{figure}

 \subsection{Orders of convergence: quarter of annulus domain}
\begin{figure}
\centering
\begin{subfigure}{0.4\textwidth}
    \includegraphics[height=8cm,width=8cm]{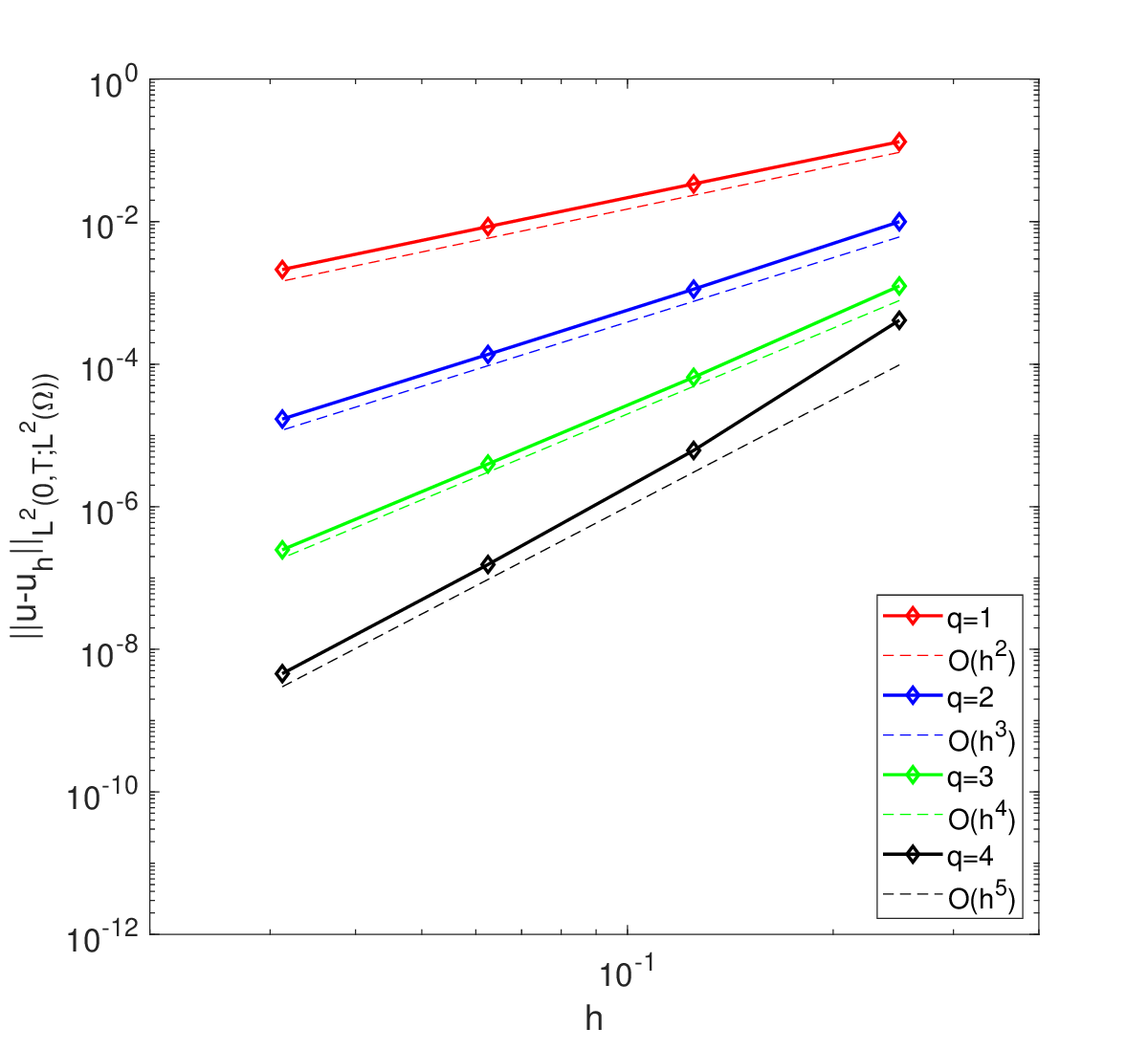}
    
    \caption{$L^2(0,T;L^2(\Omega))$ norm errors.}
    \label{fig:ringL2}
\end{subfigure}
\hspace{1cm}
\begin{subfigure}{0.4\textwidth}
    \includegraphics[height=8cm,width=8cm]{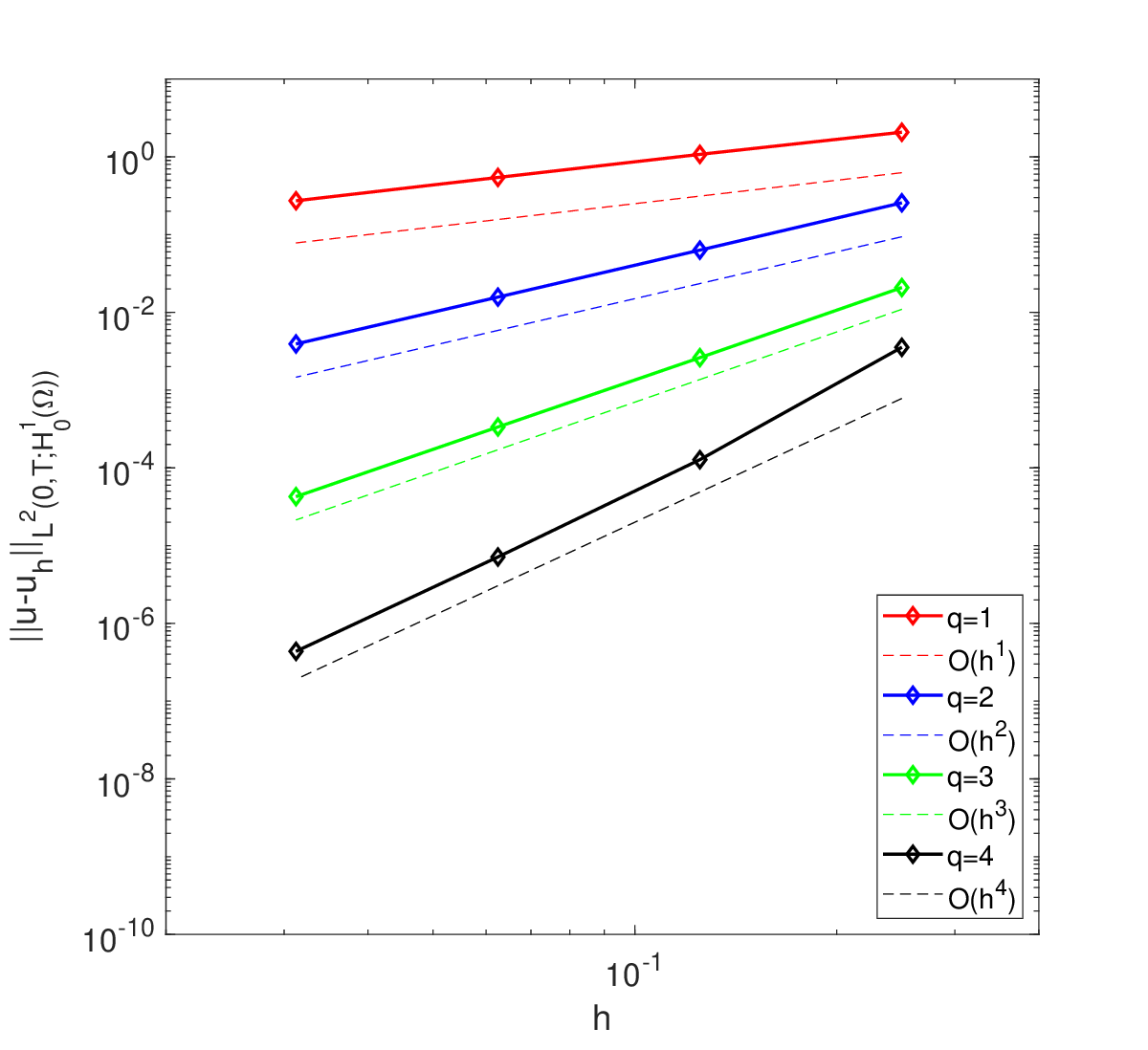}

    \caption{$L^2(0,T;H_0^1(\Omega))$ norm errors.}
    \label{fig:ringH1}
\end{subfigure}
\caption{Errors in $L^2(0,T;L^2(\Omega))$ and $L^2(0,T;H^1_0(\Omega))$ norm for quarter of annulus domain.}
\label{fig3.2}
\end{figure}

 We consider as spatial domain $\Omega$ the quarter of annulus with interior radius equal to 1 and exterior radius equal to 2 (see Figure \ref{fig:firstring}). The nonlocal term for this example is taken to be $\displaystyle a(l(u))=2-\frac{1}{1+(l(u))^2}$. We set the exact solution $u(x,t)=(x^2+y^2-1)(x^2+y^2-4)xy\sin(t)$  and compute the corresponding $f$ from \eqref{eqn1.1}. We consider basis functions of degrees $q=1,\ 2,\ 3,\ 4$ on meshes with size $h=\frac{1}{4},\ \frac{1}{8},\ \frac{1}{16}, \ \frac{1}{32}$. In Figure \ref{fig3.2}, the errors in $L^2(0,T;L^2(\Omega))$ and $L^2(0,T;H_0^1(\Omega))$ norms are plotted with respect to $h$. From Figure \ref{fig3.2}, it can be observed that the order of convergence in $L^2(0,T;H_0^1(\Omega))$ is $O(h^q)$ which coincides with the theoretical estimates given in the Theorem \ref{thm6.1}. We also compute the order of convergence in $L^2(0,T;L^2(\Omega))$ norm, even if not proved in our theorems, and we find that it is of order $O(h^{q+1})$ which is still optimal.		

 \subsection{Orders of convergence: thick ring-shaped domain}\label{secthick_ring}
\begin{figure}
\centering
\begin{subfigure}{0.4\textwidth}
    \includegraphics[height=8cm,width=8cm]{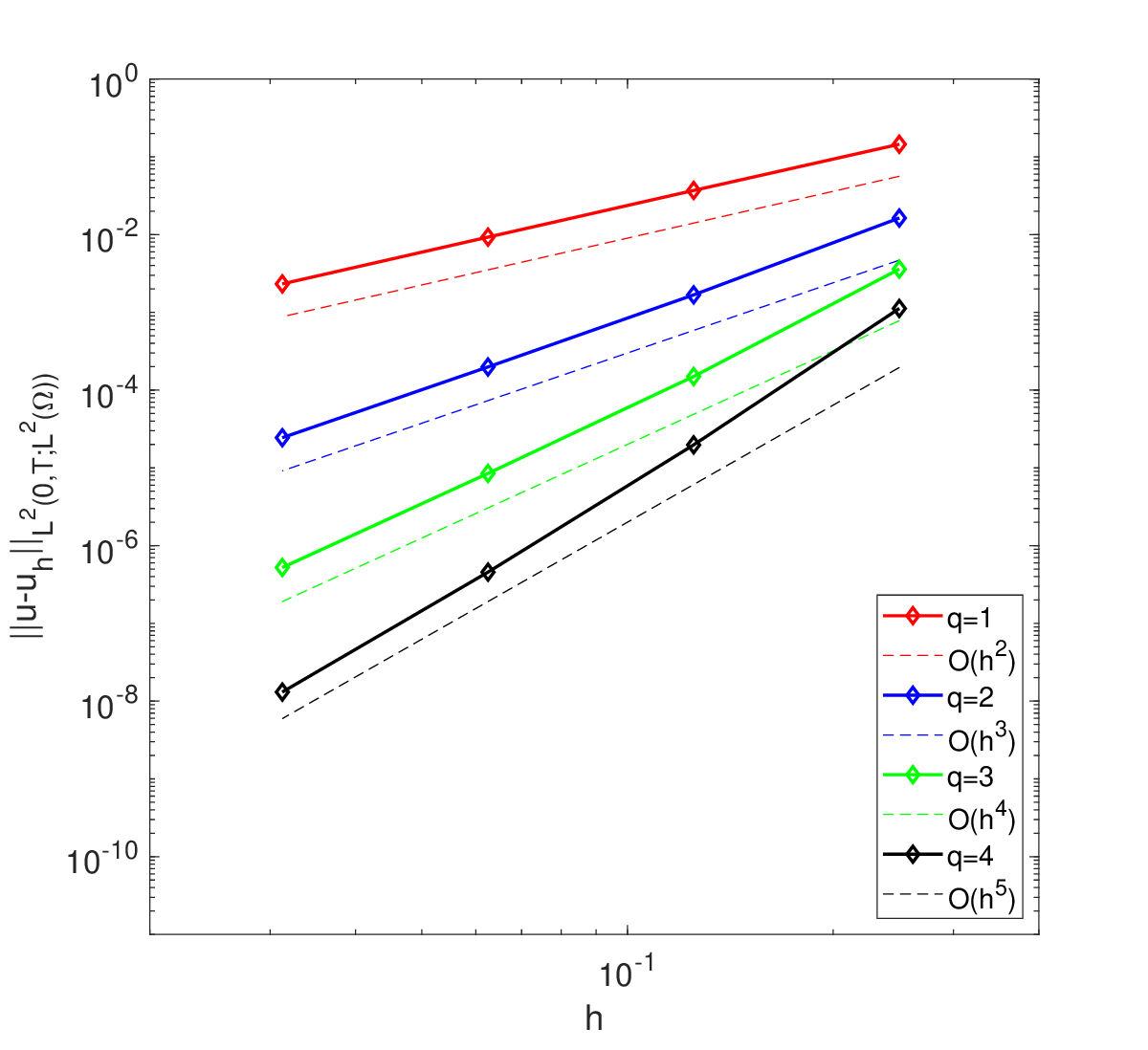}
    
    \caption{$L^2(0,T;L^2(\Omega))$ norm errors.}
    \label{fig:thick_ringL2}
\end{subfigure}
\hspace{1cm}
\begin{subfigure}{0.4\textwidth}
    \includegraphics[height=8cm,width=8cm]{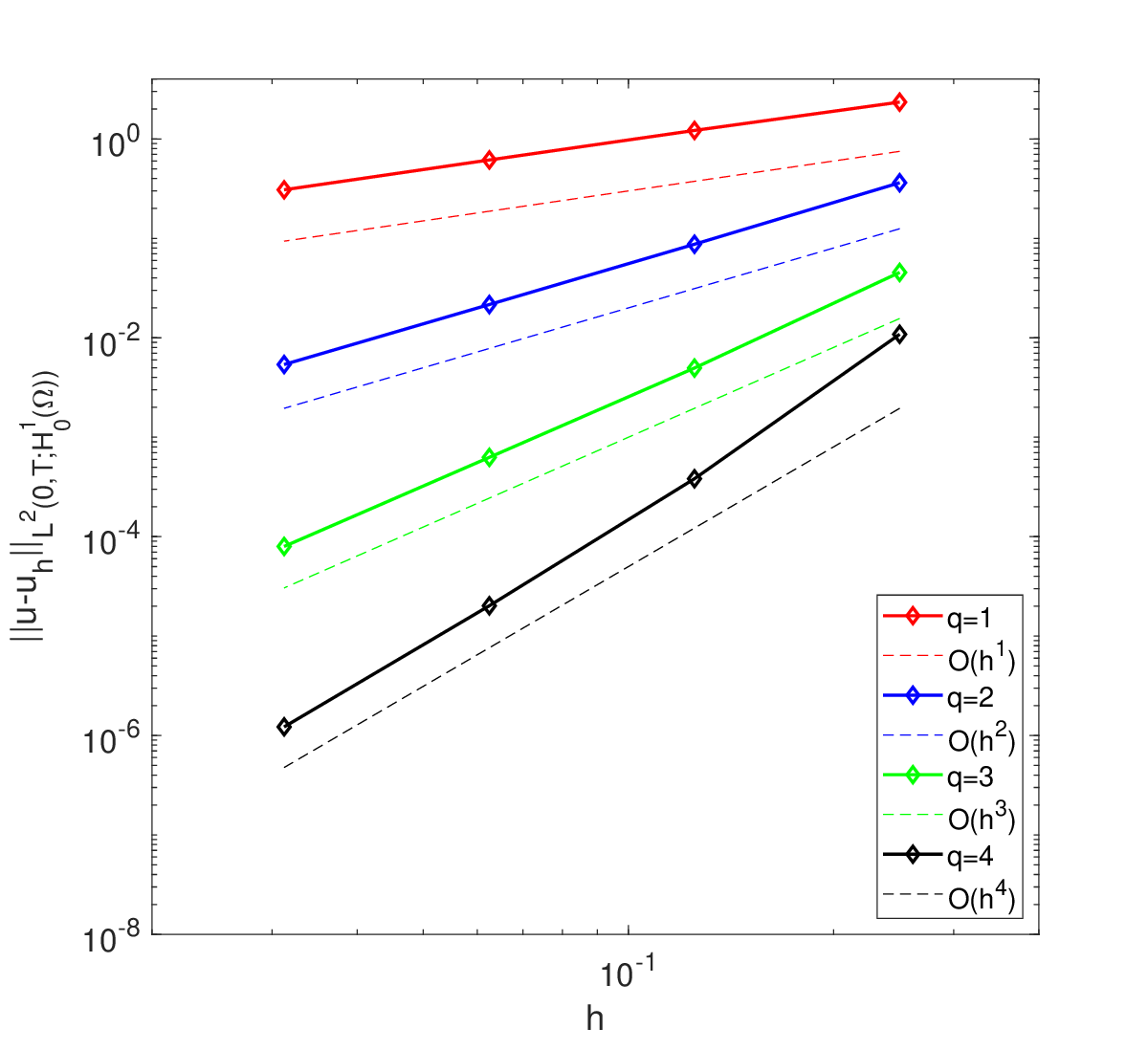}

    \caption{$L^2(0,T;H_0^1(\Omega))$ norm errors.}
    \label{fig:thick_ringH1}
\end{subfigure}
\caption{Errors in $L^2(0,T;L^2(\Omega))$ and $L^2(0,T;H^1_0(\Omega))$ norm for the thick ring-shaped  domain.}
\label{fig:thick_ringOC}
\end{figure}

In this example, we consider as computational domain the thick ring-shaped domain with interior radius equal to 1, exterior radius equal to 2 and height equal to 1, as in Figure  \ref{fig:thick_ring}. The right-hand function $f$ is taken so that the exact solution is $u(x,y,z,t)=-(x^2+y^2-1)(x^2+y^2-4)\sin(t)\sin(\pi z)xy^2$ for the nonlocal term  $a(l(u))=3+\sin(l(u))$. We consider basis functions of degrees $q=1,\ 2,\ 3,\ 4$ on meshes with  size $h=\frac{1}{4},\ \frac{1}{8},\ \frac{1}{16}, \ \frac{1}{32}$. In Figure \ref{fig:thick_ringOC}, we report the convergence rates in $L^2(0,T;L^2(\Omega))$ and $L^2(0,T;H_0^1(\Omega))$ norms for different degrees of B-spline basis functions. We see that the order of convergence in $L^2(0,T;L^2(\Omega))$ norm is $O(h^{q+1})$ and the order of convergence in $L^2(0,T;H_0^1(\Omega))$ norm is $O(h^q)$.
\subsection{Comparison of time-stepping and space-time method}
\begin{figure}[H]
		\centering
		\begin{subfigure}{0.4\textwidth}
			\includegraphics[height=8cm,width=8cm]{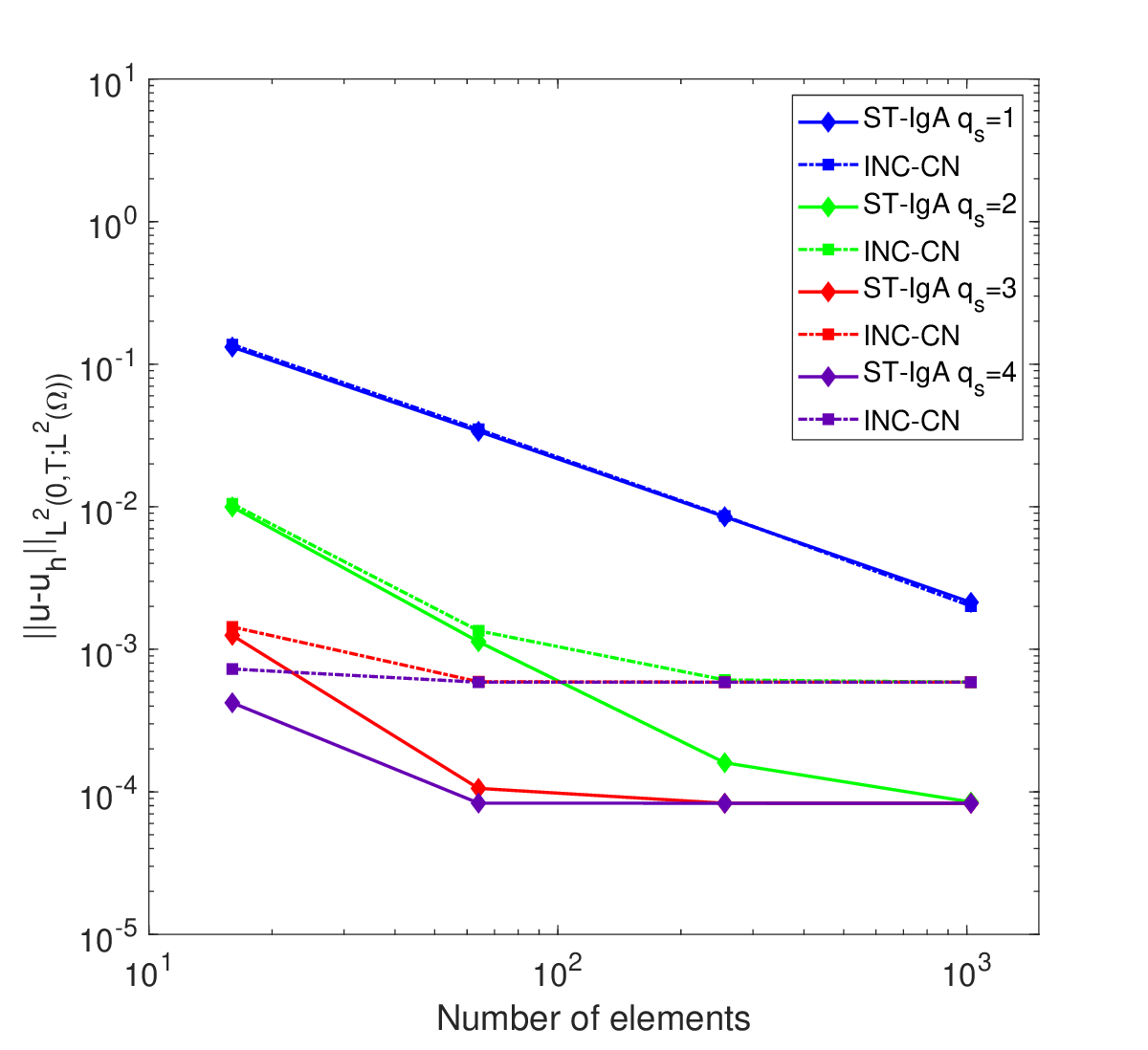}
			
			\caption{$L^2(0,T;L^2(\Omega))$ norm errors with $h_t=\displaystyle\frac{1}{16}$ and $q_t=1$.}
			\label{}
		\end{subfigure}
		\hspace{1cm}
		\begin{subfigure}{0.4\textwidth}
			\includegraphics[height=8cm,width=8cm]{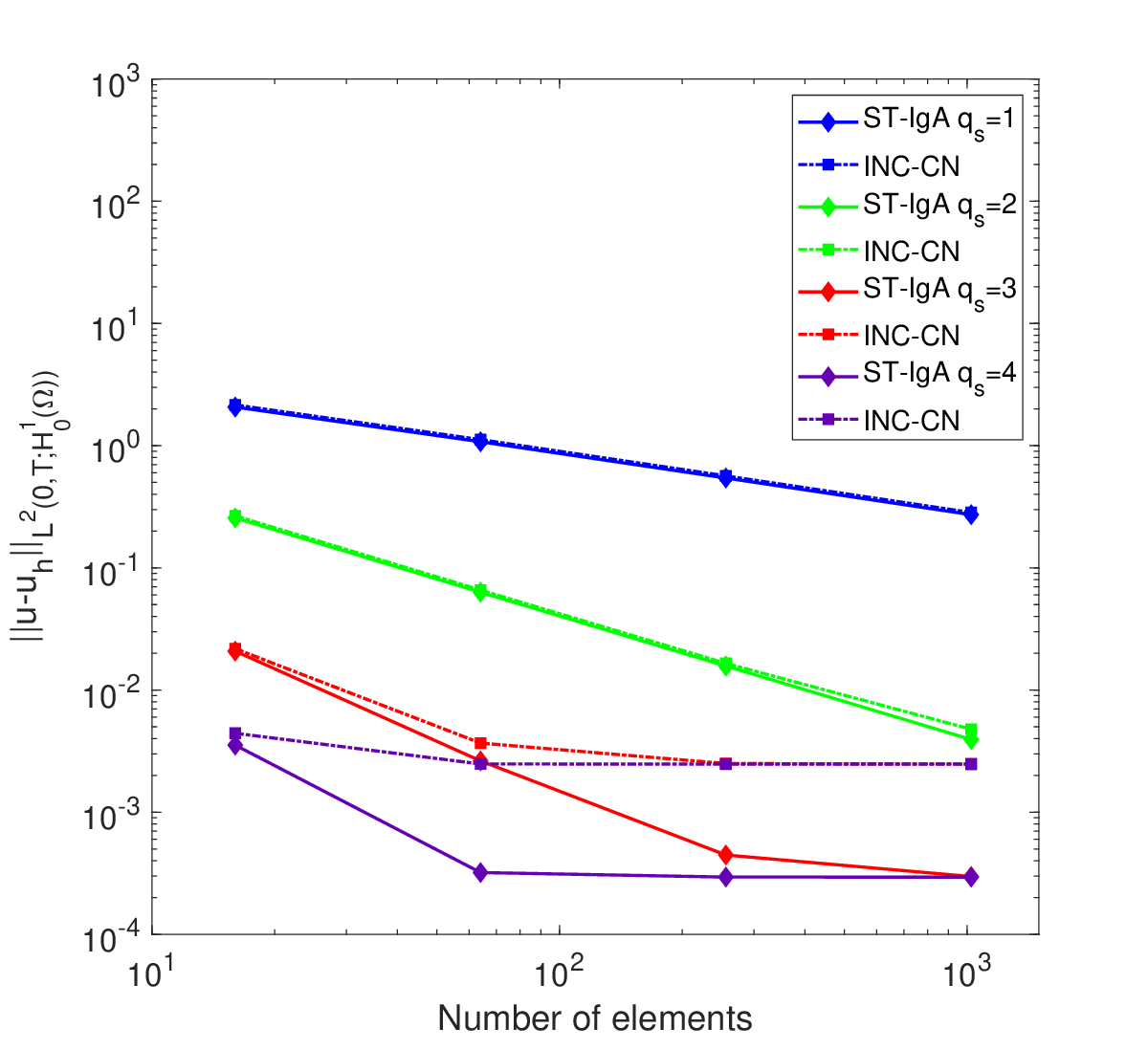}
			
			\caption{$L^2(0,T;H_0^1(\Omega)) $ norm errors with $h_t=\displaystyle\frac{1}{16}$ and $q_t=1$.}
			\label{}
		\end{subfigure}
		\caption{Errors in $L^2(0,T;L^2(\Omega))$ and $L^2(0,T;H^1_0(\Omega))$ norm for INC-CN and ST-IgA approaches}
		\label{comp1}
\end{figure}

\begin{figure}[H]
	\centering
	\begin{subfigure}{0.4\textwidth}
		\includegraphics[height=8cm,width=8cm]{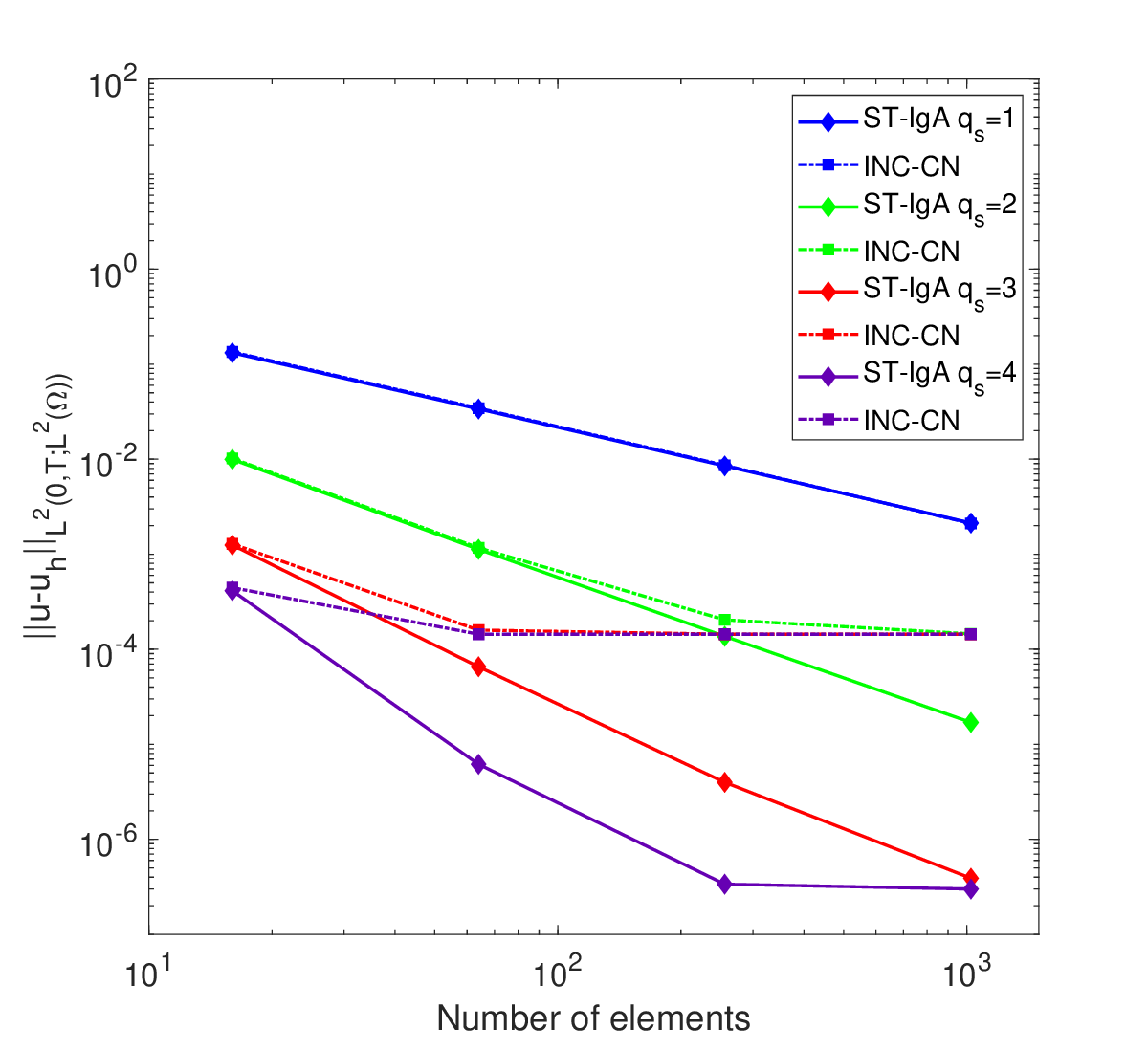}
		
		\caption{$L^2(0,T;L^2(\Omega))$ norm errors with $h_t=\displaystyle\frac{1}{32}$ and $q_t=2$.}
		\label{}
	\end{subfigure}
	\hspace{1cm}
	\begin{subfigure}{0.4\textwidth}
		\includegraphics[height=8cm,width=8cm]{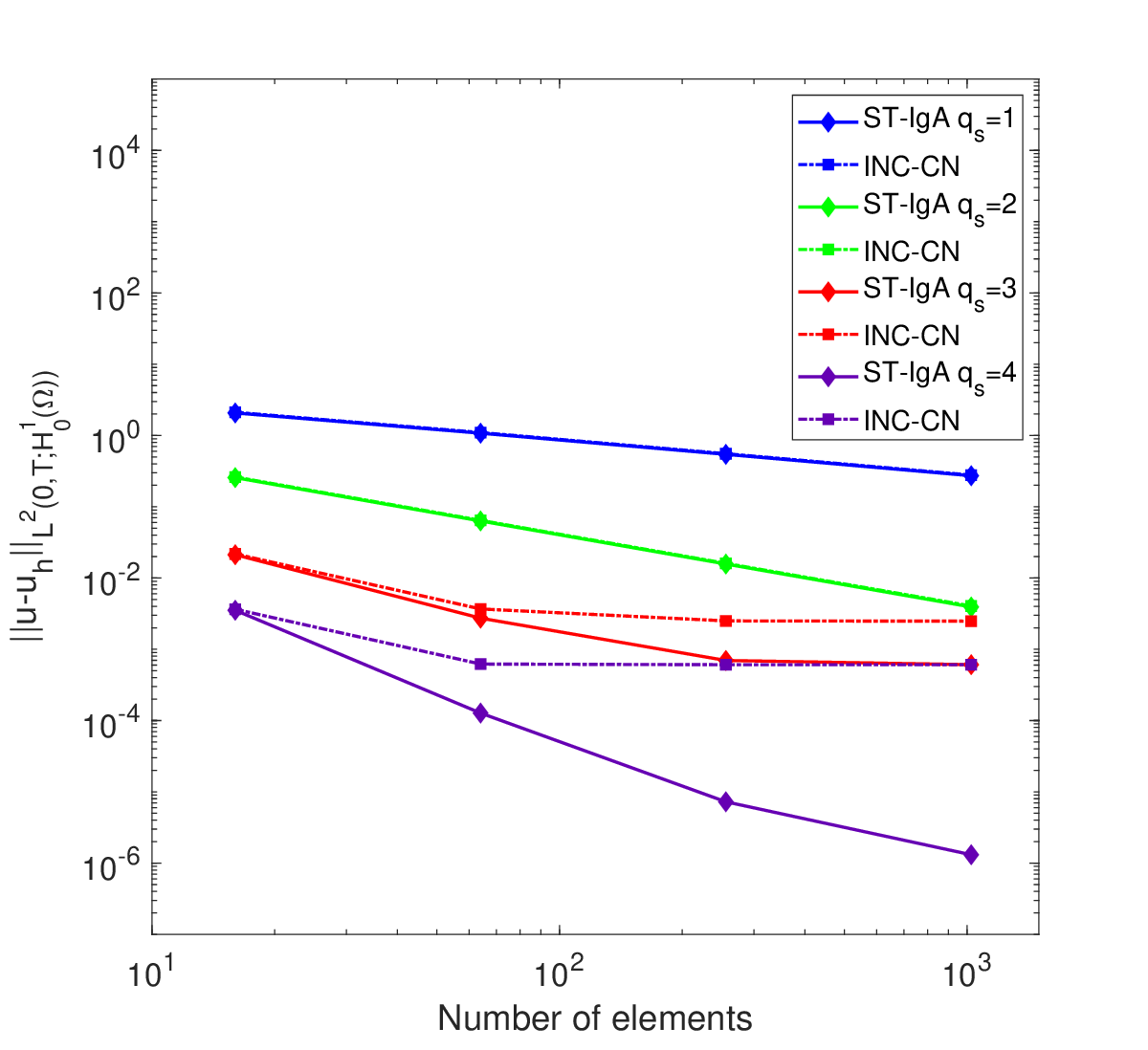}
		
		\caption{$L^2(0,T;H_0^1(\Omega))$ norm errors with $h_t=\displaystyle\frac{1}{32}$ and $q_t=2$.}
		\label{}
	\end{subfigure}
	\caption{Errors in $L^2(0,T;L^2(\Omega))$ and $L^2(0,T;H^1_0(\Omega))$ norm for INC-CN and ST-IgA approaches}
	\label{comp2}
\end{figure}

\begin{figure}[H]
	\centering
	\includegraphics[height=10cm,width=10cm]{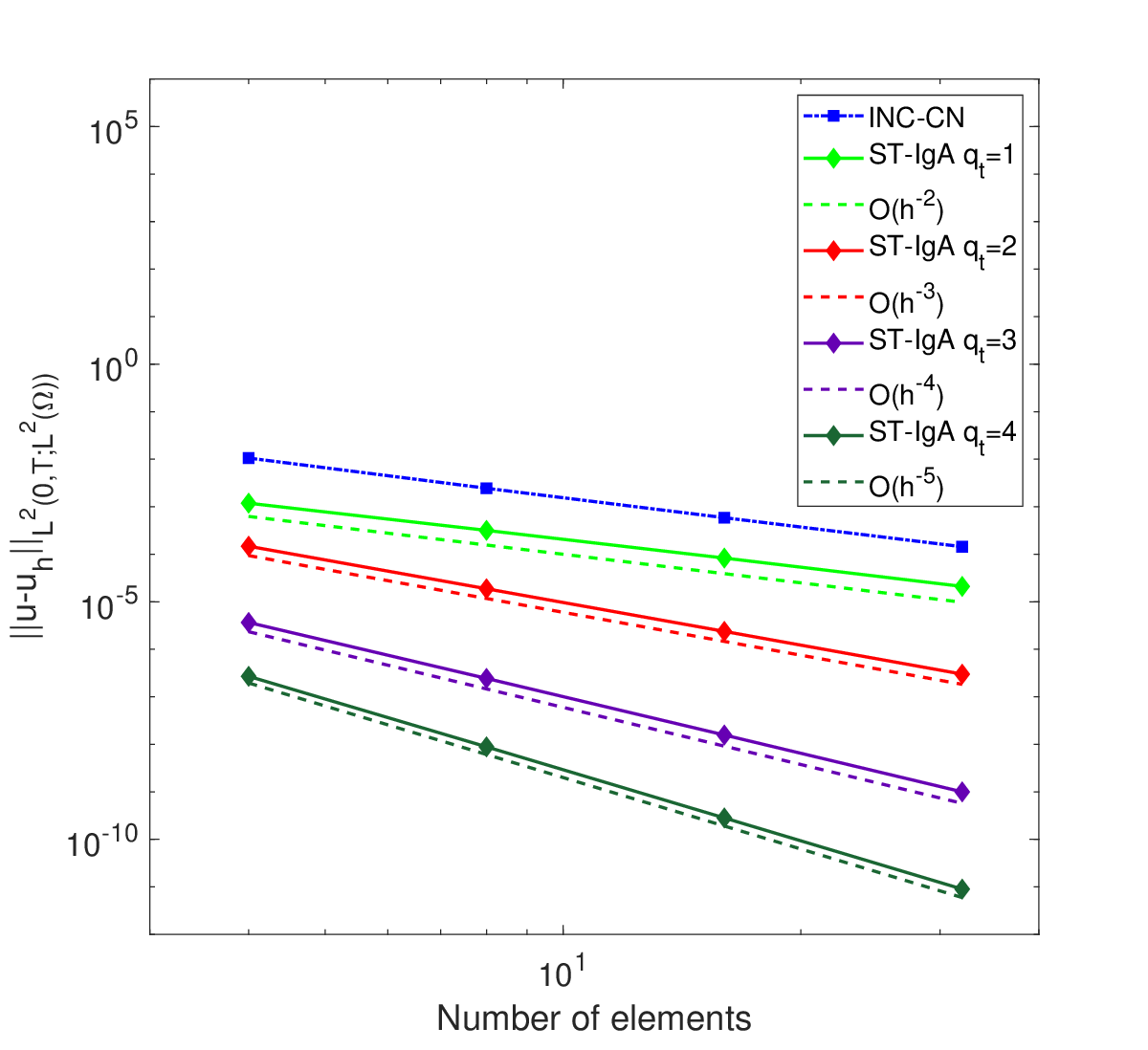}
	\caption{ Errors in $L^2(0,T;L^2(\Omega))$ norm for INC-CN and ST-IgA approaches}
	\label{cmp3}
\end{figure}

	We present an example to compare the proposed space-time isogeometric approach with a traditional time-stepping scheme, i.e. we use Crank-Nicolson method to handle time discretization. A similar comparison has also been done in \cite{nonlinearstiga} for nonlinear time-dependent problems. \\
    We consider as spatial domain $\Omega$ the quarter of an annulus with interior radius equal to 1 and exterior radius equal to 2 (see Figure \ref{fig:firstring}). The nonlocal term for this example is taken to be $\displaystyle a(l(u))=2-\frac{1}{1+(l(u))^2}$. We set the exact solution $u(x,y,t)=(x^2+y^2-1)(x^2+y^2-4)xy\sin(t)$  and compute the corresponding $f$. \\    
     We solve the nonlinear problem using both the Crank-Nicolson incremental approach \cite{INC-CN} (denoted by INC-CN) and our proposed space-time isogeometric method (denoted by ST-IgA). Specifically, we perform our numerical experiments for two different cases: (i) we fix the temporal meshsize $h_t$ and we calculate the errors for different spatial meshsizes $h_s$ and different spatial degrees $q_s$ (ii) we fix the spatial meshsize and calculate the errors for different temporal meshsizes $h_t$ with different temporal degrees $q_t$.\\\\ 
\noindent{\bf{Case 1:}} (a) First, we set the time step size $h_t=\displaystyle\frac{1}{16}$ for the Crank-Nicolson scheme, while for the space-time approach, we fix the temporal meshsize $h_t=\displaystyle\frac{1}{16}$ and the degree $q_t=1$. We display the errors in $L^2(0,T;L^2(\Omega))$ and $L^2(0,T;H_0^1(\Omega))$ norms in Figure \ref{comp1}. We use square markers to display errors for the INC-CN approach and diamond markers to display errors for the ST-IgA approach.  In Figure \ref{comp1}, we observe that the errors obtained using the  ST-IgA approach are smaller than those of the INC-CN approach.
\hspace{-0.7cm}(b) Similarly, we also solve the nonlinear matrix system using INC-CN approach with $h_t=\displaystyle\frac{1}{32}$ while for the ST-IgA approach we set $h_t=\displaystyle\frac{1}{32}$ and $q_t=2$ and we display the errors for different spatial meshsizes in Figure \ref{comp2}. From Figure \ref{comp2}, we observe that, despite refinement, with INC-CN, we are unable to reduce the error below a certain threshold, while in the ST-IgA approach, the errors reduce as long as the meshsize reduces. This shows the superiority of the space-time approach.
\hspace{-0.6cm}\noindent{\bf{Case 2:}} In this test, we calculate the errors for different time refinements and for fixed spatial meshsize $h_s=\displaystyle\frac{1}{64}$ and degree $q_s=8$. In Figure \ref{cmp3}, we display the $L^2(0,T;L^2(\Omega))$ errors obtained using INC-CN approach for different $h_t$ and the ST-IgA approach for different meshsizes $h_t$ and different temporal degrees $q_t$. From Figure \ref{cmp3}, we observe an order of convergence $\mathcal{O}(h_t^2)$  for  INC-CN  an order of convergence of $\mathcal{O}(h_t^{q_t+1})$, $q_t=1,2,3,4$ for ST-IgA. The results show that the ST-IgA approach yields more accurate solutions than the INC-CN approach.

\subsection{Performance of preconditioner: igloo-shaped domain }
\label{sec:6.4}
In this example, we consider an igloo-shaped domain as represented in Figure \ref{fig:second}.
We set the right-hand side as  $f=1$ and the nonlocal term as $a(l(u))=3+\sin(l(u))$. We consider zero initial condition and homogeneous Dirichlet boundary conditions. Note that, in this example the exact solution is unknown and the domain $\Omega$ is not trivial. In this example, we are interested in studying the number of Picard's and the GMRES iterations which are required in finding the numerical solution. We  consider basis functions of degrees $q=1,\ 2,\ 3,\ 4$ on meshes with  size $h= \frac{1}{8},\ \frac{1}{16}, \ \frac{1}{32}, \ \frac{1}{64}$. Table \ref{tabann} shows the number of Picard's iterations and the maximum number of GMRES iterations at any one of the Picard's iterations. We also have noticed that we get the same number of GMRES iterations (up to 1 more/ 1 less) in every Picard's iteration. We can see that the number of GMRES iterations is almost independent of $q$ and $h$.
\begin{table}[H]
    \centering
    \begin{tabular}{ccccc}
    \hline
        \multicolumn{5}{c}{Picard's iterations/ GMRES iterations} \\
          \hline
       $h^{-1}$\hspace{1cm} & $q=1$\hspace{1cm} & $q=2$\hspace{1cm}  & $q=3$\hspace{1cm} & $q=4$\\
       \hline
        8\hspace{1cm} & 3/14\hspace{1cm} & 3/14\hspace{1cm} & 3/15\hspace{1cm} & 3/16 \\
         16\hspace{1cm} & 3/18\hspace{1cm} & 3/18\hspace{1cm} & 3/19\hspace{1cm} & 3/19 \\
         32\hspace{1cm} & 3/20\hspace{1cm} & 3/20\hspace{1cm} & 3/20\hspace{1cm} & 3/21 \\
         64\hspace{1cm} & 3/21\hspace{1cm} & 3/21\hspace{1cm} & 3/22\hspace{1cm} &  3/22 \\
         \hline
    \end{tabular}
    \caption{Igloo-shaped domain performance of preconditioner}
    \label{tabann}
\end{table}
\subsection{A numerical experiment with different nonlocal term}
\begin{figure}[H]
	\centering
	\begin{subfigure}{0.4\textwidth}
		\includegraphics[height=8cm,width=8cm]{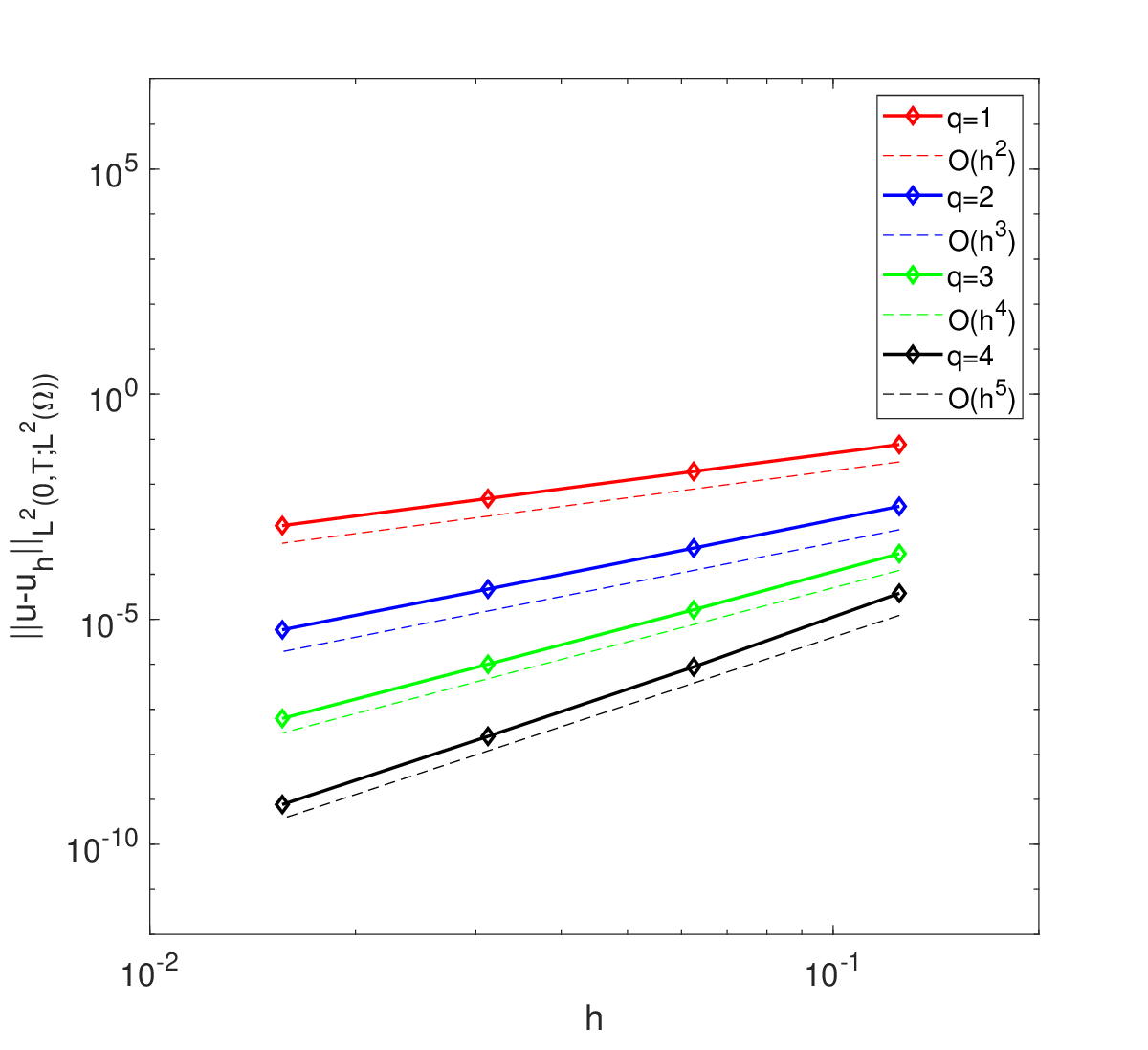}
		
		\caption{$L^2(0,T;L^2(\Omega))$ norm errors.}
		\label{fig:ring_compL2}
	\end{subfigure}
	\hspace{1cm}
	\begin{subfigure}{0.4\textwidth}
		\includegraphics[height=8cm,width=8cm]{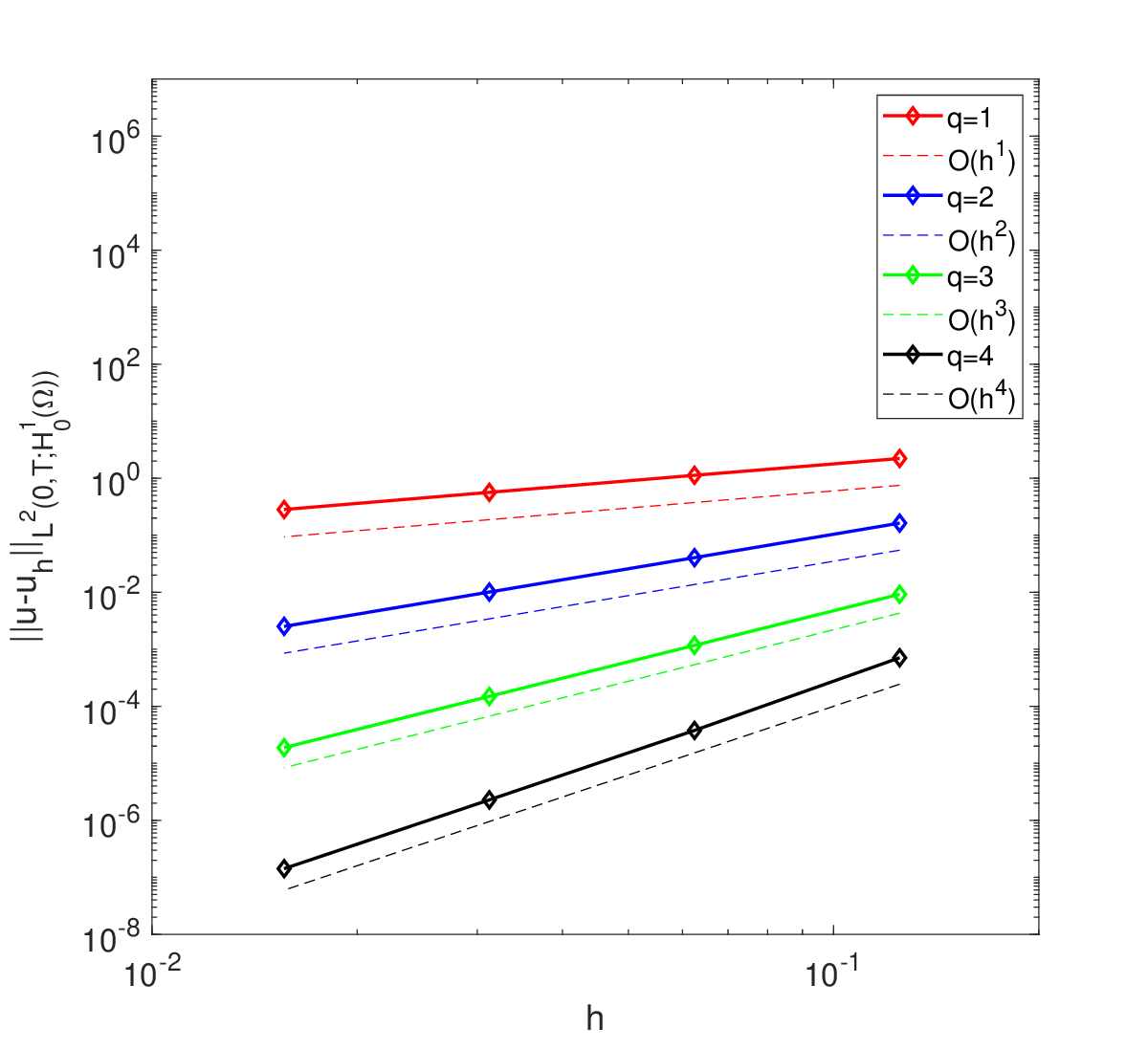}
		
		\caption{$L^2(0,T;H_0^1(\Omega))$ norm errors.}
		\label{fig:ring_compH1}
	\end{subfigure}
	\caption{Errors in $L^2(0,T;L^2(\Omega))$ and $L^2(0,T;H^1_0(\Omega))$ norm for quarter of annulus domain.}
	\label{diffnonlocal}
\end{figure}

In this example, we consider problem \eqref{eqn1.1} with a different nonlocal term $l(u)$. In particular, we choose $l(u)=\displaystyle\int_{\Omega}\vert\nabla u\vert^2\ d\Omega$  and we solve   numerically the problem using the proposed space-time isogeometric method, even if  this type of non-locality  is not covered in our theoretical analysis. We choose as domain $\Omega$ a quarter of annulus with interior radius equal to 1 and exterior radius equal to 2 (see Figure \ref{fig:firstring}) and take $\displaystyle a(l(u))=2-\frac{1}{1+(l(u))^2}$. The right-hand side function $f$ is calculated so that the exact solution is $u(x,y,t)=-xy^2\sin(\pi t)(x^2 + y^2 - 1)(x^2 + y^2 - 4)$. In Figure \ref{diffnonlocal}, the errors in $L^2(0,T;L^2(\Omega))$ and $L^2(0,T;H_0^1(\Omega))$ norms are displayed for meshsizes $h=\frac{1}{8}, \frac{1}{16}, \frac{1}{32},\frac{1}{64}$ and polynomial degrees $q=1,2,3,4$. We observe that the order of convergence in $L^2(0,T;L^2(\Omega))$ norm is $\mathcal{O}(h^{q+1})$ and  in $L^2(0,T;H_0^1(\Omega))$ norm is $\mathcal{O}(h^q)$. In Table \ref{tabring}, we display the number of Picard's iterations and the maximum number of GMRES iterations at any one of the Picard's iterations. We can see that the number of GMRES iterations is almost independent of $q$ and $h$. As in the test case of Section \ref{sec:6.4}, we observed that we get the same number of GMRES
iterations (up to 1 more/ 1 less) in every Picard’s iteration.
\begin{table}[H]
    \centering
    \begin{tabular}{ccccc}
    \hline
        \multicolumn{5}{c}{Picard's iterations/ GMRES iterations} \\
          \hline
       $h^{-1}$\hspace{1cm} & $q=1$\hspace{1cm} & $q=2$\hspace{1cm}  & $q=3$\hspace{1cm} & $q=4$\\
       \hline
        8\hspace{1cm} & 30/11\hspace{1cm} & 13/11\hspace{1cm} & 14/12\hspace{1cm} & 19/12 \\
         16\hspace{1cm} & 12/12\hspace{1cm} & 19/12\hspace{1cm} & 20/12\hspace{1cm} & 16/13 \\
         32\hspace{1cm} & 22/12\hspace{1cm} & 16/12\hspace{1cm} & 16/13\hspace{1cm} & 16/13 \\
         64\hspace{1cm} & 16/12\hspace{1cm} & 15/12\hspace{1cm} & 15/13\hspace{1cm} &  15/13 \\
         \hline
    \end{tabular}
    \caption{Quarter of annulus domain performance of preconditioner}
    \label{tabring}
\end{table}

\section{Conclusions}
In this work, we have proposed space-time isogeometric method for a parabolic problem with nonlocal diffusion coefficient. We have shown the existence-uniqueness of the solution at continuous level. The existence of a discrete solution is proved using the Brouwer's fixed point theorem. However, proof of the uniqueness of the discrete solution is yet to be discovered. Theoretically, we have derived the error estimate in $L^2(0, T; H_0^1(\Omega))$ norm and we have confirmed it through the numerical experiments. Although the error estimate in $L^2(0, T; L^2(\Omega))$ norm is not proven theoretically, the optimal order of convergence in $L^2(0, T; L^2(\Omega))$ norm is observed in the numerical experiments.
We have also provided a preconditioner suited for the linearized space-time isogeometric discretization of the nonlocal problem, based on \cite{LOLI20202586}. At each non-linear iteration, the linear system matrix $\mathbf{A}(\mathbf{u})$ changes and the preconditioner must be built from scratch at each iteration. However, the set-up and application cost of preconditioner is almost optimal as discussed in Section \ref{sec4}. Numerical experiments have confirmed that the proposed preconditioner is robust with respect to the degree and the meshsize.
Although our benchmark examples do not show signs of numerical instability, unphysical behaviors and spurious temporal oscillations can occur.  A detailed investigation of these issues is left for future works.\\\\

{\noindent{\bf CRediT authorship contribution statement:}}\\\\ \textbf{Sudhakar Chaudhary:} Conceptualization, Formal analysis, Methodology, Validation, Supervision, Writing-original draft, \textbf{Shreya Chauhan:} Conceptualization, Formal analysis, Methodology, Software, Validation, Writing- original draft, \textbf{Monica Montardini:} Formal analysis, Methodology, Software, Validation, Supervision, Writing- original draft \\\\\\
{\noindent{\bf Acknowledgment:}} The authors acknowledge anonymous reviewers for many helpful suggestions and comments. The second author would like to express gratitude to the Department of Science and Technology (DST), New Delhi, India for providing INSPIRE Fellowship to carry out the research work. The third author is member of the Gruppo Nazionale Calcolo
Scientifico-Istituto Nazionale di Alta Matematica (GNCS-INDAM).\\\\

{\noindent {\bf Data availability:} }The data that support the findings of this study are available from the corresponding author upon reasonable request.\\\\
{\noindent {\bf Declarations}}\\\\
{\bf Conflict of interest:} On behalf of all authors, the corresponding author states that there is no conflict of interest.\\\\


\end{document}